\documentclass[reqno,final]{amsart}
\usepackage{natbib}  
\usepackage{fancyhdr} 
\usepackage{color} 
\usepackage{hyperref} 
\usepackage{graphicx} 

\definecolor{aleacolor}{rgb}{0.16,0.59,0.78}

\hypersetup{
breaklinks,
colorlinks=true,
linkcolor=aleacolor,
urlcolor=aleacolor,
citecolor=aleacolor}

\renewcommand{\headheight}{11pt}
\renewcommand{\cite}{\citet}

\theoremstyle{plain}
\newtheorem{thm}{Theorem}[section]  
\newtheorem{prop}[thm]{Proposition}   
\newtheorem{lem}[thm]{Lemma} 
\newtheorem{cor}[thm]{Corollary}  

\theoremstyle{definition}
\newtheorem{defn}[thm]{Definition} 
\theoremstyle{remark}
\newtheorem{rem}[thm]{Remark} 
\newtheorem{ex}[thm]{Example} 

\makeatletter \@addtoreset{equation}{section} \makeatother
\renewcommand\theequation{\thesection.\arabic{equation}}
\renewcommand\thefigure{\thesection.\arabic{figure}}
\renewcommand\thetable{\thesection.\arabic{table}}

\newcommand{\aleaIndex}[1]{\href{http://alea.impa.br/english/index_v#1.htm}{\bf #1}}
\eheader{Alea}{\aleaIndex{7}}{2010}{79}{116}

\elogo{\parbox[c]{3cm}{\includegraphics[width=3cm]{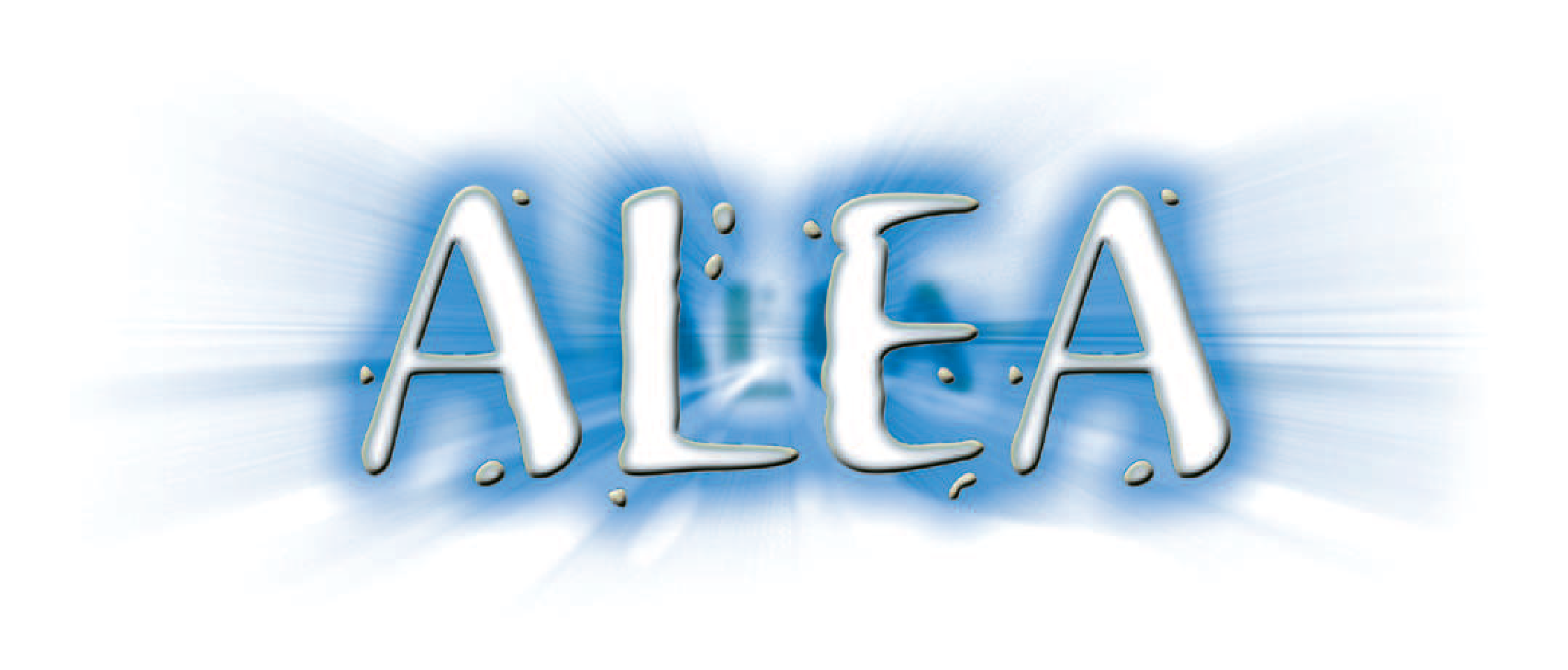}}}

\usepackage{natbib}
\usepackage{fancyhdr}
\usepackage{graphicx}
\usepackage{amssymb}
\usepackage{amsmath}

\newcommand{\cal}{\mathcal}
\newcommand{\al}{\alpha}

\newcommand{\gm}{\gamma}

\newcommand{\ep}{\varepsilon}

\newcommand{\sg}{\sigma}

\newcommand{\ps}{\psi}

\newcommand{\Om}{\Omega}

\newcommand{\mcal}{\mathcal}

\newcommand{\q}{\quad}
\newcommand{\qq}{\qquad}

\newcommand{\la}{\langle}
\newcommand{\ra}{\rangle}

\newcommand{\R}{\mathbb{R}}
\newcommand{\N}{\mathbb{N}}
\newcommand{\Z}{\mathbb{Z}}

\newcommand{\les}{\leqslant}
\newcommand{\ges}{\geqslant}

\def\prf{\noindent{\textit{Proof: }}}

\begin{document}

\date{April 24, 2009; accepted February 23, 2010}
\keywords{selfsimilar process, selfdecomposability, dominated variation, submultiplicativity, O-subexponentiality} 
\subjclass{60G18, 60G51, 60E07.}

\author{Toshiro Watanabe}
\address{Center for Mathematical Sciences, The University of Aizu, Aizu-Wakamatsu, \\
Fukushima 965-8580, Japan.}
\email{t-watanb@u-aizu.ac.jp}

\author{Kouji Yamamuro}
\address{Faculty of Engineering, Gifu University, Gifu 501-1193, Japan.}
\email{yamamuro@gifu-u.ac.jp}

\title[]{Limsup behaviors of multi-dimensional selfsimilar processes with independent increments} 

\begin{abstract}
Laws of the iterated logarithm of ``limsup'' type are studied for multi-dimensional selfsimilar  processes $\{X(t)\}$ with independent increments  having exponent $H$. It is proved that, for any positive increasing function $g(t)$ with $\displaystyle \lim_{t\to\infty}g(t) = \infty$, there is $C\in [0,\infty]$ such that $\limsup|X(t)|/(t^Hg(|\log t|))= C$ a.s. as $t \to\infty $, in addition, as $t \to 0$. A necessary and sufficient condition for the existence of  $g(t)$ with $C=1$ is obtained. In the case where $g(t)$ with $C=1$ does not exist, a criterion to classify functions $g(t)$ according to $C=0$  or $C=\infty$ is given. Moreover, various ``limsup'' type laws  with identification of  the positive constants $C$ are explicitly presented in several propositions and examples. The problems that exchange the roles of $\{X(t)\}$ and  $g(t)$ are also discussed.
\end{abstract}

\maketitle

\section{Introduction}

Since  Lamperti (1962),  the study of selfsimilar processes has been made by many scholars.  In a celebrated paper \citet{sa91}, K. Sato introduced and characterized selfsimilar processes with independent increments without assuming the stationarity of increments. As far as the authors know, M. Yor and D. Madan, in the second L\'evy conference at Aarhus in 2002, called those processes {\it Sato processes}. This naming is  used in applications to mathematical finance, for example, in \citet{cgmy05,cgmy07}. K. Sato called those processes {\it selfsimilar additive processes} in his book \citet{sa99} afterwards. In this paper, we discuss ``limsup'' type limit theorems for multi-dimensional selfsimilar additive processes as time goes to infinity, in addition, to zero. The four problems are proposed later on. They are motivated and started from a pioneer paper \citet{pr90} for increasing random walks. Our aim is to solve them. It is known by \citet{sa91} that the marginal distribution at any time of an $H$-selfsimilar additive process is selfdecomposable and conversely, for any selfdecomposable distribution $\mu$ and $H>0$, there is a unique in law $H$-selfsimilar additive process such that the marginal distribution at time 1 is the same as the distribution $\mu$. Thus the investigation on the tail behaviors of selfdecomposable distributions is crucial for the resolution of the problems. In particular, the dominated variation and the evolution of  generalized moments of an infinitely divisible distribution play key roles.  The greatest difficulty lies in the fact that an infinitely divisible distribution and its L\'evy measure do not always have the same tail behaviors.

Historically, \citet{sa91} started to study the rate of growth of selfsimilar additive processes in some increasing case. Following this, \citet{wa96}  investigated  ``limsup'' and ``liminf'' type limit theorems for the general increasing case and \citet{ya03}  treated the same problem in a certain class of two-sided processes. As closely related results, the sample path behaviors of selfsimilar Markov processes are investigated in \citet{cp06,  pa06,  vi03}.  A selfsimilar Markov process was introduced and characterized by \citet{la72}. \citet{sy98,sy00} and \citet{ya00a,ya00}  studied recurrence-transience for selfsimilar additive processes in detail. As a remarkable fact, \cite{ya00}  proved that all selfsimilar additive processes in dimensions greater than or equal to three are transient. In dimensions $1$ and $2$, the attempt to find a criterion of their recurrence-transience has not been successful. In recent years, problems analogous to the present paper  were discussed in a series of works \cite{wa02,wa02a,wa04,wa07}  for shift selfsimilar additive random sequences. They were applied to solve some classical  problems of random fractals.

In what follows, let $\R^d$ be the $d$-dimensional Euclidean space and $S^{d-1}$ be the $(d-1)$-dimensional unit sphere with the understanding that $S^0:= \{-1,1\}$.
The symbol $|x| $ stands for the Euclidean norm of $x$ in $\R^d$ and $\langle z,x \rangle$ does for the Euclidean inner product of $z$ and $x$ in $\R^d$.  Let $\R_+ = [0,\infty)$. We use the words ``increase'' and ``decrease'' in the wide sense allowing flatness. Thus the word ``monotone'' means either ``increase'' or ``decrease''. A precise definition of a selfsimilar additive process is given as below.

\begin{defn}
  An $\R^d$-valued stochastic  process $\{X(t):t\geq 0\}$ on a probability space $(\Omega, {\mathcal F}, P)$  is called a {\it selfsimilar additive process} (or Sato process) with exponent $H>0$ if the following conditions are satisfied:

{\rm (d.1)} $\{X(t)\}$ is selfsimilar with exponent $H$, that is,
$$\{X(ct)\}\stackrel{\rm d}{=} \{c^HX(t)\}\mbox{ for every $c>0$,}$$
where the symbol  $\stackrel{\rm d}{ =}$ stands for the equality of finite dimensional distributions.

{\rm (d.2)} $\{X(t)\}$ has independent increments, that is, 
$$X(t_1)-X(t_0),\ X(t_2)-X(t_1),\cdots,X(t_{n})-X(t_{n-1})$$
 are independent for any $n$ and any choice of $0\leq t_0< t_1< \cdots< t_n$.

{\rm (d.3)} $\{X(t)\}$ has c\`adl\`ag path, that is,  almost surely $X(t)$ is right-continuous in $t\geq 0$ and has left limits in $t>0$.\par
Instead of a selfsimilar additive process with exponent $H$, we sometimes say an {\it $H$-selfsimilar additive process}. 
\end{defn}

Note that the stochastic continuity of  $\{X(t)\}$ follows from (d.1) and (d.2) above and that we do not assume the stationarity of increments. If the process  $\{X(t)\}$ has  stationary increments, that is, if the process is a L\'evy process, then it is a strictly $\alpha$-stable L\'evy process with $\alpha=1/H$. If there is $c\in \R^d $ such that $X(t) = t^Hc$, then $\{X(t)\}$ is called deterministic. Otherwise it is called non-deterministic. The process $\{X(t)\}$ is called symmetric if $\{X(t)\}\stackrel{\rm d}{=} \{-X(t)\}$. The process  $\{X(t)\}$ is called Gaussian if the distribution of $X(t)$ is Gaussian at any time $t>0$, otherwise it is called non-Gaussian. 

Before \citet{sa91}, \citet{ge79}  gave  the following two examples of selfsimilar additive processes in relation to the $d$-dimensional Brownian motion.  Let $\{B(t), t \geq 0\}$ be  a Brownian motion  on  $\R^d$, starting at the origin. Define the hitting time process  $\{T_r, r  \geq 0 \}$ and, for $d \ge 3$, the last exit time process  $\{L_r, r  \geq 0 \}$  as follows:  For  $r \ge 0$, 
$$ T_r: = \inf \{t \ge 0 : |B(t)| = r \} \mbox{ and  }  L_r: = \sup \{t \ge 0 : |B(t)|= r \}.$$ 
Then it is proved by \citet{ge79} that   $\{T_r, r  \geq 0 \}$ and $\{L_r, r  \geq 0 \}$ are $\R_+$-valued increasing $2$-selfsimilar additive processes.  Note that in the case $d=3$, $\{L_r, r  \geq 0 \}$ is a one-sided 1/2-stable L\'evy process and, except this case, $\{T_r, r  \geq 0 \}$ and $\{L_r, r  \geq 0 \}$ are not L\'evy processes. The explicit representations of the distributions of $T_1$ and $L_1$ are already known by \citet{ct} and \citet{ge79}, respectively. Here we define a class ${\mathcal G}_1$ of functions $g$ (the so-called $g$-functions) on $\R_+$  as
\begin{eqnarray}
\label{eq:1.1}
&& {\mathcal G}_1:=\{g(x): \mbox{$g(x)$ is positive and increasing on $[0,\infty)$,}\nonumber\\
&&\qquad\quad \lim_{x \to \infty}g(x) = \infty \}.
\end{eqnarray}
Directly applying the key lemmas (Lemmas 3.1, 3.2, 4.3, and 4.5) of \cite{wa96}, we obtain the following results both as $r \to 0$ and as $r \to \infty$. Assertion (i) below is also found in \cite{gs95}.

{\bf Theorem A.}\ \begin{it}  {\rm (i)} We have 
 $$ \liminf\frac{\log|\log r|}{r^2}T_r = \frac{1}{2} \quad \mbox{ a.s.} $$

{\rm (ii)} Let  $k:= d/2 -1$. We have 
 $$ \limsup\frac{T_r}{r^2\log|\log r|} = \frac{2}{j_k^2} \quad \mbox{ a.s.} $$
 where $j_k$ is the first positive zero of the Bessel function $J_k(x)$ of the first kind.

{\rm (iii)}  Let $d \ge 3$. We have 
$$ \liminf\frac{\log|\log r|}{r^2}L_r = \frac{1}{2} \quad \mbox{ a.s.} $$

{\rm (iv)}  Let $d \ge 3$ and $g \in {\mathcal G}_1$. When we consider the case $r \to 0$, we make an additional assumption that  $x^2g(|\log x|)$ is increasing in $x$ with $0<x<1$.
If $\int_0^{\infty}(g(x))^{-(d/2 -1)}dx < \infty$ (resp. $ = \infty)$, then
$$ \limsup\frac{L_r}{r^2g(|\log r|)} =  0 \mbox{ (resp. $  =\infty)$}\mbox{ a.s. }$$
\end{it}

It is obvious that the laws in Theorem A  are equivalent to the well known laws of the  iterated logarithm for the Brownian motion $\{B(t)\}$. Thus we can give another nice proof for the following  classical laws both as $t \to 0$ and as $t \to \infty$.

 {\bf Corollary A.}\ \begin{it}
 {\rm (i)}  We have 
$$ \limsup\frac{1}{\sqrt{t\log|\log t|}}\sup_{0\le s \le t}|B(s)| = \sqrt{2} \quad \mbox{ a.s.} $$

{\rm (ii)} Let  $k:= d/2 -1$. We have 
 $$ \liminf\sqrt{\frac{\log|\log t|}{t}}\sup_{0\le s \le t}|B(s)| = \frac{j_k}{\sqrt{2}} \quad \mbox{ a.s.} $$

{\rm (iii)}   Let $d \ge 3$. We have 
$$ \limsup\frac{1}{\sqrt{t\log|\log t|}}\inf_{t\le s }|B(s)| = \sqrt{2} \quad \mbox{ a.s.} $$

{\rm (iv)}  Let $d \ge 3$ and $g \in {\mathcal G}_1$. When we consider the case $t \to 0$, we make an additional assumption that  $x^2g(|\log x|)$ is increasing in $x$ with $0<x<1$.
If $\int_0^{\infty}(g(x))^{-(d/2 -1)}dx < \infty$ (resp. $ = \infty)$, then
$$ \liminf\sqrt{\frac{g(|\log t|)}{t}}\inf_{t\le s }|B(s)| =  \infty \mbox{ (resp. $ =0)$}\mbox{ a.s. }$$
\end{it}

In Corollary A,  ``$\sup_{0\le s \le t}|B(s)|$'' in (i) and ``$\inf_{t\le s }|B(s)|$'' in (iv) can be replaced by $|B(t)|$. Assertion (i), (ii), (iii) and (iv) above are due to \cite{kh33}, \cite{ta67}, \cite{kl94} and \cite{de50}, respectively. In particular, (ii) is called the Chung type law of the iterated logarithm. The order of the assertions in Theorem A corresponds to that in Corollary A. Throughout this paper, let $\{X(t)\}$ be an $H$-selfsimilar additive process on $\R^d$. We investigate limsup behaviors of $\{X(t)\}$. \par

By virtue of our key theorem (Theorem \ref{th2.1} in Sect.2), we obtain that for any  $g \in {\mathcal G}_1$, there is $C\in [0,\infty]$  such that
\begin{eqnarray}
\label{eq:1.2}
\limsup_{t \to \infty} \frac{|X(t)|}{t^Hg(\log t)} = C \quad\mbox{a.s. } 
\end{eqnarray}
When $0<C<\infty$, we take $Cg(x)$ in place of $g(x)$. Hence the cases to consider are three cases where $C=0$, $C=1$ and $C=\infty$. The law \eqref{eq:1.2} is called to be normal if $C=1$. Now we propose four problems. As is seen in (i) of Corollary A, the Gaussian case is already known. In the present paper, we focus on the non-Gaussian case.  \par
If we suppose that $\limsup_{x \to \infty} g(x+1)/g(x) < e^{H}$ for $g \in {\mathcal G}_1$, our key theorem (Theorem \ref{th2.1}) remains true with \eqref{eq:1.2} replaced by 
\begin{eqnarray}
\label{eq:1.3}
\limsup_{t\to 0}\frac{|X_t|}{t^Hg(|\log t|))}=C\quad a.s.
\end{eqnarray}
Hence our theorems and propositions hold with \eqref{eq:1.2} replaced
by \eqref{eq:1.3}.
\medskip

{\bf Problem 1.} (Normalizability)

Suppose that $\{X(t)\}$ is given. What is a necessary and sufficient condition for the existence of $g \in {\mathcal G}_1$ satisfying \eqref{eq:1.2} with $C=1$ ?
Further, in what way can we give this  $g \in {\mathcal G}_1$ ?
\medskip
\goodbreak

{\bf Problem 2.} (Integral test for $g$)

Suppose that $\{X(t)\}$ is given and that  there is no $g \in {\mathcal G}_1$ satisfying \eqref{eq:1.2} with $C=1$.   What criterion classifies   functions $g \in {\mathcal G}_1$ into those satisfying  \eqref{eq:1.2}  with $C=0$ and  those satisfying  \eqref{eq:1.2} with $C=\infty$ ?
\medskip

{\bf Problem 3.} (Possible types of $g$-functions in the normal  laws) 

Suppose that $g \in {\mathcal G}_1$  is given. 
What is a necessary and sufficient condition in terms of  $ g \in {\mathcal G}_1$ for the existence of $ \{X(t)\}$   satisfying \eqref{eq:1.2} with $C=1$ ?
Moreover,  in what way can we give this $ \{X(t)\}$  ?
\medskip

{\bf Problem 4.} (Integral test for $ \{X(t)\}$)

Suppose that $g \in {\mathcal G}_1$  is given and that  there is no  $ \{X(t)\}$ satisfying  \eqref{eq:1.2} with $C=1$. What criterion classifies processes  $ \{X(t)\}$  into those satisfying  \eqref{eq:1.2}  with $C=0$ and  those satisfying  \eqref{eq:1.2} with $C=\infty$ ?
\medskip

In Sect.2, we answer Problems 1 and 2 in Theorems \ref{th2.2} and
\ref{th2.3}, and partially but substantially answer Problems 3 and 4
in Theorems \ref{th2.4}, \ref{th2.5} and \ref{th2.6}. In addition, we
present several interesting and explicit examples with identification
of the constant $C$ in \eqref{eq:1.2}. In Sect.3, we show some
preliminary results on the tail behaviors of multivariate infinitely
divisible distributions. In Sect.4, we prove the main results
mentioned in Sect.2. In Sect.5, we prove Theorem \ref{th3.1} in
Sect.3.


\medskip

\section{Answers to four problems}

 In this section, we answer the problems mentioned in Sect.1. We give the answers only as $t \to \infty$. Those to the case $t \to 0 $ are similar and omitted.  \par
 
Let $\N =\{1,2, \ldots \}$, $\Z=\{0,\pm 1,\pm 2, \ldots \}$, and $\Z_+=\{0, 1, 2, \ldots \}$. 
 For positive measurable functions $f(x)$ and $g(x)$ on $\R^1$, we define the relation $f(x) \sim g(x)$ by $\lim_{x \to \infty}f(x)/g(x) =1$ and the relation $f(x) \asymp g(x)$ by $0 <\liminf_{x \to \infty}f(x)/g(x) \le\limsup_{x \to \infty}f(x)/g(x) < \infty$. We denote the tail of a measure $\eta$ on $\R^d$ by $\eta(|x| >r)$, that is,  $\eta(|x| >r):= \eta(\{x:|x| >r\})$ and the right tail of a measure $\zeta $ on $\R_+$ by $\zeta(x >r)$, namely,  $\zeta(x >r):= \zeta(\{x:x >r\})$.  For a measure $\eta$ on $\R^d$, we define the probability measure $\widetilde\eta$  by 
\begin{equation}
\label{eq:2.1}
\widetilde\eta(dx) :=(\eta(|x| >1))^{-1}1_{\{|x| >1\}} \eta(dx)
\end{equation}
 only when $0 <\eta(|x| >1)< \infty$. Denote by $\delta_a(dx)$ the Dirac mass at $a \in \R^d$, that is, the probability measure concentrated at  $a \in \R^d$.  Let $\mu$ and $\rho$ be distributions on $\R^d$. Denote by $\mu*\rho$ the convolution of  $\mu$ and $\rho$ and by $\mu^{n*}$ the $n$th convolution power of $\mu$. Furthermore, we denote by $\widehat \mu(z)$ the characteristic function of  $\mu$. A distribution $\mu$ on $\R^d$ is called non-degenerate if its support is not included in any $(d-1)$-dimensional hyperplane in  $\R^d$. 

In what follows, we use the terminology in \citet{sa99}.  
Let  $\mu$  be an infinitely divisible distribution on $\R^d$, $d\ges 1$, with generating triplet $(A,\nu,\gm)$. Here 
$A$ is the Gaussian
covariance matrix, $\nu$ is the L\'evy measure, and $\gm$ is the location parameter.
That is, 
\begin{equation}\label{eq:2.2}
 \widehat \mu(z):= \int_{\R^d}\exp (i\la z, x\ra)\mu(dx) = \exp(\ps(z)),\qq z\in\R^d
\end{equation}
with
\begin{equation}\label{eq:2.3}
\ps(z)=\int_{\R^d}(e^{i\la z,x\ra}-1-1_{\{|x|\les1\}}(x) i\la z,x\ra)\nu(dx)
+i\la\gm,z\ra -\frac12 \la Az,z\ra,
\end{equation}
where $\nu$ is a measure on $\R^d$ satisfying $\nu(\{0\})=0$ and 
$\int_{\R^d}(1\land |x|^2)\nu(dx)<\infty$, $\gm\in\R^d$, and $A$ is a 
nonnegative-definite matrix. If $\nu \ne 0$, then $\mu$ is said to be  non-Gaussian. Further, if $A=0$, then $\mu$ is said to be purely non-Gaussian. 

\begin{defn}
A distribution $\mu$ on $\R^d$ is said to be {\it selfdecomposable} (or of class $L$) if, for each $b \in (0,1)$, there is a distribution $\lambda_b$ on $\R^d$ such that
\begin{equation}\label{eq:2.4}
\widehat \mu(z)=\widehat \mu(bz)\widehat \lambda_b(z).
\end{equation}
\end{defn}

Note that $\lambda_b$ in \eqref{eq:2.4} is infinitely divisible. Selfdecomposable distribution were introduced by \citet{le}. They are infinitely divisible and their convolutions are selfdecomposable again. Stable including Gaussian, Pareto, log-normal, logistic, gamma, F, t, hyperbolic, half-Cauchy distributions are  known to be selfdecomposable. In addition, so is a Weibull distribution with parameter $0<\alpha\leq 1$ ( see Remark 8.12 of \citealp{sa99}, as to parameter $\alpha$). The proofs of the selfdecomposability of those distributions are not trivial, because the L\'evy measures are not always explicitly known. See Example 15.13 and Exercise 34.14 of \citet{sa99} and Examples 2.4 and 9.16 in Chapter V and 12.8 in Chapter VI of \citet{sh04}. In the following lemma, assertions (i) and (ii) is due to \citet{wo80} for $d=1$ and \citet{sa80} for  $d\ge2$, and assertion (iii) is due to \citet{sa82}. 

\begin{lem}
\label{lem2.1}
 Let $\mu$ be a selfdecomposable distribution on $\R^d$.

{\rm (i)}  Let $\mu$ be non-Gaussian but not purely non-Gaussian. Then  there are Gaussian distribution  $\mu_1$ and purely non-Gaussian selfdecomposable distribution  $\mu_2$ such that $\mu = \mu_1*\mu_2$.

{\rm (ii)} Let $\mu$ be purely non-Gaussian. Then the L\'evy measure $\nu$ is expressed as
\begin{eqnarray}
\label{eq:2.5}
\nu(B)=\int_{S^{d-1}}\sg(d\xi)\int_0^\infty 1_{B}(r\xi)k_\xi(r)r^{-1}dr
\end{eqnarray}
for a Borel set $B$ in $\R^d$. Here $\sg$ is a finite non-zero  measure on $S^{d-1}$ and $k_\xi(r)$ is a nonnegative function  which is  measurable in $\xi\in S^{d-1}$ and decreasing in $r>0$. In the $\alpha$-stable case with $0<\alpha <2$, $k_\xi(r) =  r^{-\alpha}$.

{\rm (iii)}  If $\mu$ is non-degenerate, then it is absolutely continuous with respect to the Lebesgue measure.  

\end{lem}

The selfsimilar additive processes $ \{X(t)\}$ were characterized by \cite{sa91} as follows.

\begin{lem}
\label{lem2.2}
{\rm (i)} Let $\{X(t)\}$ be an $H$-selfsimilar additive process on $\R^d$. For any $t>0$, the distribution of $X(t)$ is selfdecomposable.

{\rm (ii)} Fix $H>0$. For any selfdecomposable  distribution $\mu$, there is a unique in law $H$-selfsimlar additive process $\{X(t)\}$ with distribution $\mu$ at time $1$.
\end{lem}

In the rest of this paper, we assume that $\{X(t)\}$ is non-deterministic. Denote by $\mu$ the distribution of $\{X(1)\}$. In the non-Gaussian case, the L\'evy measure $\nu$ of $\mu$ is expressed as \eqref{eq:2.5}.  The process  $\{X(t)\}$ is called non-degenerate if  $\mu$ is non-degenerate. For $l \in \N$, we define
\begin{eqnarray}
\label{eq:2.6}
\rho_l(dx):=P(X(1)-X(e^{-l})\in dx).
\end{eqnarray} 
Then $\rho_l$ is infinitely divisible. The equation \eqref{eq:2.4} holds with $b=e^{-lH}$ and $\lambda_b= \rho_l$. Let $\eta_l$  be the L\'evy measure of $\rho_l$. Note that $\int_{|x|>1} \log |x|\rho_l(dx)<\infty$, that is, $\rho_l$ has finite log-moment by Theorem 25.3 of \cite{sa99}. If $\mu$ is non-degenerate, then it is  absolutely continuous by Lemma \ref{lem2.1} (iii). Then we denote by $p(x)$ the density of $\mu$.  Let $m(d\xi)$ be the uniform probability measure on $S^{d-1}$. Let $l\in \N$. Define four functions $G_l(r)$, $K(r)$, $L(r)$, and  $M(r)$ on $(0,\infty)$ as follows:

\begin{eqnarray}
\label{eq:2.7}
G_l(r):=P(|X(1)-X(e^{-l})|>r)=\rho_l(|x| >r),
\end{eqnarray} 

\begin{eqnarray}
\label{eq:2.8}
K(r):=\int_{S^{d-1}}k_\xi(r)\sg(d\xi),
\end{eqnarray}

\begin{eqnarray}
\label{eq:2.9}
L(r):= P(r<|X(1)|\le e^{3H}r) = \mu(r < |x| \le e^{3H}r),
\end{eqnarray} 
and in the case where the measure $\mu$ is non-degenerate, we define 
 \begin{eqnarray}
\label{eq:2.10}
M(r):= \int_{S^{d-1}}p(r\xi)m(d\xi).
\end{eqnarray}
Then the tail of the L\'evy measure $\eta_l$ of $\rho_l$ is expressed as
\begin{eqnarray}
\label{eq:2.11}
\eta_l(|x| >r) = \int_r^{e^{lH}r}\frac{K(s)}{s} ds.
\end{eqnarray} 
As $\{X(t)\}$ is non-deterministic, $G_l(r)$ and $L(r)$ are positive on $(0,\infty)$ but $K(r)$ can be 0 for some $r >0$. In particular, it means the Gaussian case that $K(r)=0$ for all $r>0$.

\begin{defn}
\label{def2.2}
 Let ${\mathcal K}=\R^d$ or let ${\mathcal K}=\R_+$. A positive measurable function $h(x)$ on ${\mathcal K}$ is called {\it submultiplicative} on ${\mathcal K}$ if there is  $ c>0$ such that $h(x+y)\leq ch(x)h(y)$ for  every $x,  y\in {\mathcal K}$. 
\end{defn}

\begin{rem}
\label{rem2.1}
(i)
If $h(x)$ is submultiplicative on $\R_+$, then so is   $h(|x|)$   on $\R^d$.

(ii) If $h(x)$ is submultiplicative on $\R^d$, then there are $c_1,
c_2> 0$ such that $h(x)\le c_1\exp(c_2 |x|)$ for  every $x \in
\R^d$. See Lemma 25.5 of \cite{sa99}. For example, we consider a
function $h(x)=\exp[c|x|^{\alpha}(\log( |x|+1))^{\beta}]$ on $\R^d$,
where $c>0$, $\alpha>0$ and $-\infty<\beta<\infty$. Then $h(x)$ is submultiplicative only in the case where $0<\alpha < 1$ or the case where $\alpha =1$ and $-\infty < \beta \le 0$. 
\end{rem}

\begin{defn}
\label{def2.3}
 A positive measurable function $h(x)$ on $\R_+$ is said to belong to the class $ O R$ if $h(cx)\asymp h(x)$ for any $c>0$. In particular, $h(x)$ is said of {\it dominated variation} if it is  monotone and $h(2x)\asymp h(x)$. Then we write $h\in D$.
\end{defn}

Here we remark that  $h(x)\in  D$ if and only if $h(x)$ is monotone and $h(x) \in   O R$.

\begin{defn}
\label{def2.4}
 The functions
$$ h(x) = \int_0^x q(t)dt, \q f(x)= \int_0^x q^{-1}(t)dt,$$
on $\R_+$ are called {\it Young conjugate} functions if $ q(t)$ is positive on $(0,\infty)$, right-continuous and increasing with $q(0)=0$ and $\lim_{t \to \infty} q(t)=\infty$.
\end{defn}

 For example, refer to pages 54 and 65 in \citet{bgt87} and to \citet{ka82} as to Definitions \ref{def2.3} and \ref{def2.4}, respectively. Lastly, we introduce the important classes of distributions.

\begin{defn}
\label{def2.5}
Let $\zeta$ be a distribution on $ \R_+$. 

(a) The distribution $\zeta$ is called to belong to the class ${\mathcal S}$ if $\zeta*\zeta(x >r) \sim 2\zeta(x >r)$ as $r \to \infty$.

(b) The distribution $\zeta$ is called to belong to the class ${\mathcal O}{\mathcal S}$ if $\zeta*\zeta(x >r) \asymp \zeta(x >r)$ as $r \to \infty$.
\end{defn}

\begin{rem}
A distribution $\rho$ on $ \R_+$ in the class ${\mathcal S}$ is called subexponential. Distributions in the class ${\mathcal S}$ are contained in the class ${\mathcal O}{\mathcal S}$.  See also \citet{ekm97,sw05}.
\end{rem}

\begin{defn}
\label{def2.6}
Let $\rho$ be a distribution on $ \R^d$. 

(a)  Set $\zeta(x >r):=\rho(|x| >r) $ for any $r \ge 0$. The distribution $\rho$ is called to belong to the class ${\mathcal O}{\mathcal S}$ if $\zeta\in {\mathcal O}{\mathcal S}$ as a distribution on $\R_+$.

(b)  The distribution $\rho$ is called to belong to the class ${\mathcal D}$ if $\rho (|x| >r) \in D.$
\end{defn}

Recall the definition \eqref{eq:1.1} of the class  ${\mathcal G}_1$. For $g\in{\mathcal G}_1$, we define
\begin{equation}
g^{-1}(x):=  \sup\{ y : g(y) < x \} \mbox{ for $x \ge 0$}
\end{equation}
with the understanding that $\sup \emptyset =0$. Now we present the key theorem.

\begin{thm}
\label{th2.1}
 Let $g \in {\mathcal G}_1$.

{\rm (i)} There is $ C\in [0,\infty]$ such that \eqref{eq:1.2} holds.

{\rm (ii)} The equality \eqref{eq:1.2} holds for  some $ C\in[0,\infty]$ if and only if 
\begin{eqnarray}
\int_0^{\infty} G_l(\delta g(x))dx& &=\int_{\R^d} g^{-1}(\delta^{-1}|x|)\rho_l(dx)\nonumber\\
& &\left\{\begin{array}{l}
<\infty \quad   \mbox{ for } \delta>C\ \mbox{and all $l\in {\N}$},\\
=\infty \quad   \mbox{ for } 0< \delta<C\ \mbox{and some $l=l(\delta)\in {\N}$}.

\end{array}
\label{eq:2.13}
\right.
\end{eqnarray}

{\rm (iii)} Suppose that $g \in{\mathcal G}_1$ and $g^{-1}(|x|)$ is submultiplicative on $\R^d$. Then we have  \eqref{eq:1.2} if and only if 
\begin{eqnarray}
\label{eq:2.14}
\int_1^{\infty} \frac{K(r)-K(e^Hr)}{r}g^{-1}(\delta^{-1} r)dr
 \left\{
\begin{array}{ll}
<\infty &\q \mbox{ for } \delta > C,\\
=\infty &\q \mbox{ for }  0 < \delta <C.
\end{array}
\right.
\end{eqnarray}
In the case where ${\displaystyle \limsup_{r \to \infty}K(e^H r)/K(r)
  <1}$, \eqref{eq:2.14} is also equivalent to the following :
\begin{eqnarray}
\label{eq:2.15}
\int_1^\infty K(r)\frac{g^{-1}(\delta^{-1} r)}{r}dr
 \left\{
\begin{array}{ll}
<\infty & \mbox{ for }\delta > C,\\
=\infty &\mbox{ for } 0 <  \delta <C.
\end{array}
\right.
\end{eqnarray}
\end{thm}

\begin{rem}
\label{rem2.3}
 (i) Suppose further that ${\displaystyle\limsup_{x\to\infty}g(x+1)/g(x)<e^H}$ is satisfied for $g\in{\mathcal G}_1$. By virtue of Proposition \ref{prop4.2} in Sect.4, the theorem holds with  \eqref{eq:1.2} replaced by \eqref{eq:1.3}. 

 (ii) We can choose $ g\in{\mathcal G}_1$ such that $\int_{\R^d}
 g^{-1}(|x|)\rho_1(dx)= \infty$, because the support of $\rho_1$ is
 unbounded. It follows  from Theorem \ref{th2.1} (ii) that
\begin{eqnarray}
\limsup_{t\to\infty}\frac{|X(t)|}{t^H} = \infty\quad \mbox{a.s.} 
\end{eqnarray}
If ${\displaystyle \lim_{x \to \infty}g(x)<\infty}$, then \eqref{eq:1.2} always holds with $C=\infty$.  Hence any problem is not left by assuming that ${\displaystyle \lim_{x \to \infty}g(x)=\infty}$ in the definition of ${\mathcal G}_1$.

(iii) We see from \eqref{eq:2.13} that we can replace $ g(\log t)$ in \eqref{eq:1.2} by $ g(a\log t)$ for any $a >0$ without changing the value of the constant $C$.
\end{rem}

We answer the first part of Problem 1. In the non-Gaussian case, the tail behavior of $ G_l(r)$ is determined by that of $\nu$. Thus we see from Lemmas \ref{lem2.1} and  Theorem \ref{th2.1} that the answer should be given in terms of  the pair $(\sg,k_{\xi})$ or of the function $K(r)$. Moreover, by virtue of Lemma \ref{lem2.2},  the answer could be given by means of $\mu$.  \par 
Let $F(r)$ be arbitrarily chosen out of $G_1(r)$, $K(r)$, or $L(r)$. Note that the condition $K(r)\not\in  OR$ includes the case where $K(r)$ vanishes on $(c,\infty)$ with some $c \geq 0$. In particular, the case $c=0$ is the Gaussian case.

\begin{thm}
\label{th2.2}
 There exists 
 $ g\in{\mathcal G}_1$ such that \eqref{eq:1.2} holds with $C=1$ if and only if  $F(r)\not\in  O R$.
\end{thm}


We answer Problem 2.

\begin{thm}
\label{th2.3}
Let $g\in{\mathcal G}_1$. Suppose that $F(r)\in O R$.
Then \eqref{eq:1.2} holds with $C=0$ or $C=\infty$ according as
\begin{eqnarray}
\int_0^\infty F(g(x))dx &<& \infty\quad \mbox{or}\quad =\infty.
\end{eqnarray}

\end{thm}

\medskip

As a corollary, we give the strictly stable  L\'evy case which has been already shown by \cite{kh38}  for $d=1$ and by \cite{ya05}  for $d\ge 2$. See also \cite{pt83} for  the non-strict case with $d=1$. Denote by $\log_{(n)}r$ $n$-fold iteration of the logarithmic function. Here we define a function $L_{(n,\al,\ep)}(r)$ for $n \in \N$, $\al >0$, and $\ep \in \R^1$ by
\begin{equation}
L_{(n,\al,\ep)}(r):=(\log_{(n+1)} r)^{\ep+ 1/\alpha}\prod_{k=1}^{n}(\log_{(k)} r)^{1/\alpha}.
\end{equation}

\begin{cor}
\label{cor2.1}
Let $g\in{\mathcal G}_1$ and $\{X(t)\}$ be a strictly $\alpha$-stable  L\'evy process with  $0< \alpha< 2$. Then \eqref{eq:1.2} holds with $C=0$ or $C=\infty$ according as 
\begin{eqnarray}
\int_0^\infty (g(x))^{-\alpha}dx<\infty \quad\mbox{or}\quad =\infty.
\end{eqnarray}
Here $K(r) = cr^{-\alpha}$ with some $c >0$ and $H= 1/\alpha$.

 In particular, we have for $n \in \N$
\begin{eqnarray}
 \limsup_{t\to\infty}\frac{|X(t)|}{t^{1/\alpha}L_{(n,\al,\ep)}(t)}
\left\{\begin{array}{l}
=0 \mbox{ a.s.}\quad  \mbox{ for } \ep >0, \\
=\infty \mbox{ a.s.}\quad  \mbox{ for }\ep \le0.
\end{array}
\right.
\end{eqnarray}
\end{cor}

\medskip 

By virtue of Theorems \ref{th2.2} and \ref{th2.3} and Proposition
\ref{prop3.3} in Sect.3, we obtain the following corollary. It is
useful in the case where the density of $\mu$ is explicitly known with
$d=1$.

\begin{cor} 
\label{cor2.2}
 Suppose either that $d=1$ or that $\{X(t)\}$ is symmetric and non degenerated with $d \ge 2$. 

{\rm (i)} There is  $ g\in{\mathcal G}_1$ such that \eqref{eq:1.2} holds with $C=1$
if and only if $M(r)\not\in  OR.$

{\rm (ii)} Let $g\in{\mathcal G}_1$ and $M(r)\in OR$. Then \eqref{eq:1.2} holds with $C=0$ or $C=\infty$ according as 
\begin{eqnarray}
\int_0^\infty (g(x))^d M(g(x))dx<\infty\quad\mbox{or}\quad =\infty.
\end{eqnarray}
\end{cor}

\medskip 

We give an example of Corollary \ref{cor2.2}. In the example, the density of L\'evy measure is not simple. It was given by \cite{ha79}.

\begin{ex}
\label{ex2.1}
  If $\mu$ is a $t$-distribution on $\R^1$ with parameter $ m>0$, namely,
\begin{eqnarray*}
\mu(dx)=\Gamma\left(\frac{m+1}{2}\right)\left(\sqrt{\pi}\Gamma\left(\frac{m}{2}\right)\right)^{-1}\frac{1}{(1+x^2)^{(m+1)/2}}dx. 
\end{eqnarray*}
When $m=1$, $\mu$ is the  Cauchy distribution.
Let $g \in{\mathcal G}_1$. Then \eqref{eq:1.2} holds with $C=0$ or $C=\infty$ according as $\int_0^{\infty} (g(x))^{-m}dx <\infty $ or $=\infty$. In particular, we have, for $n \in \N$,
\begin{eqnarray*}
\limsup_{t\to\infty}\frac{|X(t)|}{t^H L_{(n,m,\ep)}(t)}=\left\{
\begin{array}{ll}
0 \quad \mbox{a.s.}& \mbox{ for }\ep>0,\\
\infty \quad \mbox{a.s.}& \mbox{ for }\ep\leq 0.
\end{array}
\right.
\end{eqnarray*}
\end{ex}

\medskip

We answer Problem 3, where the answer is partial but substantial. If $g^{-1}(|x|)$ is not submultiplicative on $\R^d$, the problem is not easy and not yet completely solved. See Theorem \ref{th2.6} below.

\begin{thm}
\label{th2.4}
 Let $g \in{\mathcal G}_1$.
 If   there exists $ \{X(t)\}$ such that \eqref{eq:1.2} holds with $C=1$,
then  $g^{-1}(x)+ \log x \notin  OR$. The converse  is also true provided that $g^{-1}(|x|)$ is submultiplicative on $\R^d$.
\end{thm}

We answer Problem 4 substantially. The answer is similar to Corollary 5.2 of \cite{wa02}.

\begin{thm}
\label{th2.5}
Let $g \in{\mathcal G}_1$ such that $g^{-1}(x)+ \log x \in  OR$. Then \eqref{eq:1.2} holds with $C=0$ or $C=\infty$ according as 
\begin{eqnarray}
\int_{{\R}^d}g^{-1}(|x|)\rho_1(dx) &<& \infty
\quad \mbox{or}\quad = \infty.
\end{eqnarray}
\end{thm}

Now we show two propositions useful for answering to the second parts of Problems 1 and 3 besides Theorem \ref{th2.1}.

\begin{prop}
\label{prop2.1}
Let $ C \in[0,\infty)$ and  $g \in{\mathcal G}_1$.

{\rm (i)}  \eqref{eq:1.2} holds if 
\begin{eqnarray}
\label{eq:2.23}
&&\int_{{\R}^d}g^{-1}(\delta^{-1} |x|)\mu(dx)
 \left\{\begin{array}{ll}
 <\infty & \q \mbox{ for } \delta >C,\\
=\infty & \q \mbox{ for }0<\delta <C.  
\end{array}
\right.
\end{eqnarray}
The converse is also true provided that $\int_{{\R}^d}g^{-1}(\delta^{-1} |x|)\mu(dx) <\infty$ for some $\delta \in (0,\infty)$.

{\rm (ii)}  Suppose that $g^{-1}(|x|)$ is submultiplicative on $\R^d$. Then \eqref{eq:1.2} holds if \eqref{eq:2.15} holds. The converse is also true provided that 
\begin{eqnarray}
\label{eq:2.24}
&&\int_1^{\infty}K(r)g^{-1}(\delta^{-1} r)r^{-1}dr <\infty\quad \mbox{for some $\delta \in (0,\infty)$}.
\end{eqnarray}

\end{prop}

\begin{rem}
\label{rem2.4}
 By virtue of Proposition \ref{prop2.1}, if \eqref{eq:2.13} holds with $C=\infty$, then \eqref{eq:2.23} holds with $C=\infty$. However, even if \eqref{eq:2.13} holds for some $C$ with $0\leq C<\infty$, its $C$ does not always satisfy \eqref{eq:2.23}. For example, let $g^{-1}(x) = \log(x+1)$ for $x \ge 1$. Then we can choose $\mu$ such that \eqref{eq:2.13} holds with $C=0$ but \eqref{eq:2.23} does with $C=\infty$. Moreover, there is $K(r)$ such that $K(r) \notin OR$ but $\widetilde{\nu}\in {\mathcal D}$, that is, $\mu \in {\mathcal D}$ by Corollary \ref{cor3.1} in Sect.3. See Remark 4.1 of \cite{waya}. Thus we see from Theorem \ref{th2.2} that there is  $g \in{\mathcal G}_1$ such that \eqref{eq:2.13} holds with $0 < C< \infty$, but \eqref{eq:2.23} holds with $C=\infty$.
\end{rem}

\begin{prop}
\label{prop2.2}
Let $C \in [0,\infty]$ and  $g \in{\mathcal G}_1$.

{\rm (i)}  Suppose that $\rho_1 \in {\mathcal O}{\mathcal S}$. 
Then  \eqref{eq:1.2} holds if and only if 
\begin{eqnarray}
\label{eq:2.25}
\int_0^\infty G_1(\delta g(x))dx
\left\{\begin{array}{ll}
<\infty & \q \mbox{ for } \delta >C,\\
=\infty & \q \mbox{ for }0<\delta <C. 
\end{array}
\right.
\end{eqnarray}

{\rm (ii)} Suppose that $\widetilde\eta_1 \in {\mathcal O}{\mathcal S}$. 
 Then  \eqref{eq:1.2} holds if and only if \eqref{eq:2.14} holds.

{\rm (iii)} Suppose that $\widetilde\nu \in {\mathcal O}{\mathcal S}$. Then  \eqref{eq:1.2} holds if  \eqref{eq:2.15} holds. 
The converse is also true provided that \eqref{eq:2.24} is satisfied.
\end{prop}

Here we give an example of Propositions \ref{prop2.1} and \ref{prop2.2}.

\begin{ex}
\label{ex2.2}
  Let $C \in [0,\infty]$ and  $g \in{\mathcal G}_1$. Suppose that $\mu$ is the standard lognormal distribution,  namely,
\begin{eqnarray*}
\mu(dx) =\frac{1}{\sqrt{2\pi}x}\exp(-(\log x)^2/2)dx\mbox{ on } \R_+.
\end{eqnarray*}
Then \eqref{eq:1.2} holds if and only if \eqref{eq:2.23} holds. In particular, we have
\begin{eqnarray*}
\limsup_{t\to\infty}\frac{X(t)}{t^H\exp\left(\sqrt{ 2\log_{(2)} t}\right)} &=& 1\quad \mbox{a.s.}
\end{eqnarray*}
\end{ex}

\medskip

We investigate the law \eqref{eq:1.2} in detail in the case where the $g$-function of the law is expressed as $g(t)=(\log t)^{1/\alpha}/\varphi(\log_{(2)} t)$ with some $0<\alpha\leq 1$ and some function $\varphi$. Put $a\vee b:=\max\{a,b\}$.

\medskip

\begin{prop}
\label{prop2.3}
 Let $ C\in [0,\infty]$. Suppose that $f(x)$ is regularly varying as $ x \to \infty$ and   $\exp(r^{\al}f(\log(r\vee 1)))$  with $0 <\al \le 1$ is increasing and submultiplicative on $\R_+$.  Then
\begin{eqnarray}
\label{eq:2.26}
\qquad \limsup_{t\to\infty}\frac{|X(t)|}{t^H(\log_{(2)} t)^{1/\al}/f(\al^{-1}\log_{(3)} t)^{1/\al}} = C\quad \mbox{a.s.}
\end{eqnarray}
if and only if 
\begin{eqnarray}
\label{eq:2.27}
 \liminf_{r \to \infty}\frac{-\log K(r)}{r^{\al}f(\log r) } = C^{-\alpha},
\end{eqnarray}  
equivalently
\begin{eqnarray}
\label{eq:2.28}
 \liminf_{r \to \infty}\frac{-\log \mu(|x| >r)}{r^{\al}f(\log r) } = C^{-\alpha}.
\end{eqnarray}  
\end{prop}

\medskip
Recall the definition of Young conjugate: The functions $h(x)$ and $f(x)$ are Young conjugate if $ h(x) = \int_0^x q(t)dt$ and $f(x)= \int_0^x q^{-1}(t)dt$, where $ q(t)$ is positive on $(0,\infty)$, right-continuous and increasing with $q(0)=0$ and ${\displaystyle \lim_{t \to \infty} q(t)=\infty}$.

\medskip

\begin{prop}
\label{prop2.4}
 Let  $C \in [0,\infty)$. Suppose that $h$ and $f$ are Young conjugate functions. Then  we have
\begin{eqnarray}
\label{eq:2.29}
\limsup_{t\to\infty}\frac{|X(t)|}{t^H(\log_{(2)} t)/f^{-1}(\log_{(3)}t)}=C\quad \mbox{a.s.}
\end{eqnarray}  
if and only if
\begin{eqnarray}
\label{eq:2.30}
 \liminf_{r \to \infty}\frac{h^{-1}(-\log K(r))}{r } =  C^{-1},
\end{eqnarray}  
equivalently
\begin{eqnarray}
\label{eq:2.31}
 \liminf_{r \to \infty}\frac{-\log \rho_1(|x| >r)}{rf^{-1}(\log r) } = C^{-1}.
\end{eqnarray} 
In particular, if $h(r)$ is regularly varying with positive index, then \eqref{eq:2.30} is equivalent to
\begin{eqnarray}
\label{eq:2.32}
 \liminf_{r \to \infty}\frac{-\log \mu(|x| >r)}{rf^{-1}(\log r) } = C^{-1}.
\end{eqnarray}  
 
\end{prop}

\medskip 

\begin{prop}
\label{prop2.5}
 Let $ C\in [0,\infty]$. 
 Then we have
\begin{eqnarray}
\qquad \limsup_{t\to\infty}\frac{|X(t)|}{t^H(\log_{(2)} t)/\log_{(3)} t} = C\quad \mbox{a.s.}
\end{eqnarray}
if and only if $C =\inf\{r > 0: K(r)=0\}$ with the understanding that $\inf \emptyset = \infty$,
equivalently
\begin{eqnarray}
 \lim_{r \to \infty}\frac{-\log \mu(|x| >r)}{r\log r } = C^{-1}.
\end{eqnarray}  
\end{prop}

\medskip

\begin{rem}
\label{rem2.5}
 \ If $\nu\not= 0$, then from Proposition \ref{prop2.5} it follows that
$$\limsup_{t\to\infty}\frac{|X(t)|}{t^H(\log_{(2)} t)/\log_{(3)} t}=C_1\quad\mbox{a.s.}$$
with $ C_1\in (0,\infty]$. If $\nu=0$, then we have by the Gaussian type law of the iterated logarithm 
$$\limsup_{t\to\infty}\frac{|X(t)|}{t^H\sqrt{\log_{(2)} t}}=C_2\quad\mbox{a.s.}$$
with $C_2 \in (0,\infty)$. 
Thus there is a big difference in the ``limsup'' behaviors of $\{X(t)\}$ between the Gaussian and the non-Gaussian case.
\end{rem}

We supplement Theorem \ref{th2.4} with the following theorem. Pay attention to Remark \ref{rem2.1} (ii). In the theorem,  $g^{-1}(|x|)$ is not submultiplicative on $\R^d$ for $g \in{\mathcal G}_1$. It shows that Problem 3 is more difficult than the analogous problem in \cite{wa02}. See Theorem 5.2 of \cite{wa02}.

\begin{thm}
\label{th2.6}
 Let $g \in{\mathcal G}_1$. 

 {\rm (i)} There is   $\{X(t)\}$ such that \eqref{eq:1.2} holds with $C=1$ in the following cases :

{\rm (1)} For some $c \in (0,\infty),$ 
\begin{eqnarray}
g^{-1}(x)\asymp\exp(cx^2).
\end{eqnarray}  

{\rm (2)} For some $c \in (0,\infty),$ 
\begin{eqnarray}
g^{-1}(x)\asymp\exp(cx\log x).
\end{eqnarray}  

{\rm (3)}  Suppose that $f(x)$ is expressed as
$f(x)= \int_0^x q(t)dt$
on $[0,\infty)$, where $ q(t)$ is positive on $(0,\infty)$, right-continuous and increasing with $q(0)=0$ and ${\displaystyle \lim_{t \to \infty} q(t)=\infty}$. Further, $g(x) $ satisfies that
\begin{eqnarray}
g^{-1}(x)\asymp\exp(xf^{-1}(\log x)).
\end{eqnarray}  

{\rm (ii)} There is no $\{X(t)\}$ such that \eqref{eq:1.2} holds with $C=1$ in the following cases: 

{\rm (1)} The function $g $ satisfies that
\begin{eqnarray}
\lim_{x \to \infty}\frac{g^{-1}(x)}{\exp(cx^2)} =\infty \mbox{ for all } c\in(0,\infty). \label{Th2.6.1}
\end{eqnarray}

{\rm (2)} The function $g $ satisfies that
\begin{eqnarray}
\lim_{x \to \infty}\frac{g^{-1}(x)}{\exp(cx\log x)} =\infty\mbox{ for all  } c\in(0,\infty),\label{Th2.6.2}
\end{eqnarray} 
and
\begin{eqnarray}
\lim_{x \to \infty}\frac{g^{-1}(x)}{\exp(cx^2)} =0 \mbox{ for all } c\in(0,\infty).\label{Th2.6.3}
\end{eqnarray} 
\end{thm}

\medskip

\begin{rem}
There is not always $\{X(t)\}$ such that \eqref{eq:1.2} holds with $C=1$ even provided that $g^{-1}(x)+\log x \notin OR$. For example, consider a function $g\in{\mathcal G}_1$ such that $g^{-1}(x)\asymp\exp(c_0x(\log x)^{\al_0})$ for some $\al_0 >1$ and  $c_0\in(0,\infty)$, and apply Theorem \ref{th2.6} (ii).  
\end{rem}

Finally, we give an example of Propositions \ref{prop2.3} and \ref{prop2.4}. It is also related with Proposition \ref{prop2.5} and Theorem \ref{th2.6}. 
 The example  shows that even if $K(r)$ has the same form of functions with parameters $\al$ and $\beta$, there is a delicate difference in the law \eqref{eq:1.2} according as  the parameters change their values. The difference essentially comes from whether $g^{-1}(x)$ is submultiplicative on $\R_+$ for the $g$-function of the law or not.

\begin{ex}
\label{ex2.3}
 Let $d \ge 1$, $c, c'\in [0,\infty]$, $0< \alpha, \alpha'< \infty$, and $-\infty < \beta, \beta', \gamma' < \infty$. We consider the following two kinds of assumptions : 
\begin{eqnarray}
\label{eq:2.41}
\liminf_{r \to \infty}\frac{-\log K(r)}{r^{\alpha}(\log r)^{\beta}}=c,
\end{eqnarray}
and
\begin{eqnarray}
\label{eq:2.42}
\liminf_{r \to \infty}\frac{-\log \mu(|x| >r)}{r^{\alpha'}(\log r)^{\beta'}(\log_{(2)} r)^{\gamma'}}=c'.
\end{eqnarray}

(i) Suppose either that $0<\alpha<1$ or that $\alpha=1$ and $\beta \le 0$. Only in this case, $g^{-1}(x)$ is submultiplicative for the $g$-function of the law \eqref{eq:2.43} below. We have
\begin{eqnarray}
\label{eq:2.43}
\limsup_{t\to\infty}\frac{|X(t)|}{t^H(\log_{(2)} t)^{1/\al}/(\log_{(3)} t)^{\beta/\alpha}} &=& \alpha^{\beta/\alpha}c^{-1/\alpha}\quad \mbox{a.s.}
\end{eqnarray}
if and only if \eqref{eq:2.41} holds, equivalently \eqref{eq:2.42} holds with $c' =c$, $\alpha'=\al$,  $\beta' =\beta $, and $\gamma'= 0$.
In particular, if $\mu$ is a Weibull distribution with parameter $0 < \alpha \le 1$, namely, $\mu(dx) = \alpha x^{\alpha -1}\exp(-x^{\alpha})dx\mbox{ on } \R_+$, then \eqref{eq:2.43} holds with $c=1$ and $\beta=0$.

(ii) Suppose that $c\in (0,\infty]$ and either that $\alpha=1$ and $\beta > 0$ or that $\alpha>1$ and $-\infty < \beta < \infty$. Let $C_1:=  c^{-1}$ for $\alpha=1$ and  $C_1:= \alpha^{(\beta-\alpha)/\alpha}(\alpha-1)^{(\alpha-1)/\alpha}c^{-1/\alpha}$ for  $\alpha>1$. Then we have
\begin{eqnarray*}
\limsup_{t\to\infty}\frac{|X(t)|}{t^H(\log_{(2)} t)/((\log_{(3)} t)^{(\alpha-1)}(\log_{(4)} t)^{\beta})^{1/\alpha}} &=& C_1 \quad \mbox{a.s.}
\end{eqnarray*}
if and only if \eqref{eq:2.41}  holds, equivalently \eqref{eq:2.42} holds with $c'=C_1^{-1}$, $\alpha'=1$, $\beta' =(\alpha-1)/\alpha $, and $\gamma'= \beta/\alpha $.
\end{ex}


\section{Results on tail behaviors}

In this section, we give several preliminary results on the tail behaviors of infinitely divisible distributions on $\R^d$. The results  below except for the lemmas are new and of interest in themselves. Theorem \ref{th3.1} below is proved only in Sect.5. As in Sect.2, we denote by $\mu$ the distribution of $X(1)$. Then $\mu$ is a selfdecomposable distribution on $\R^d$ except for the delta measures. We continue to use the notation of Sect.2. The following definition is due to \cite{dj88}.

\begin{defn}
\label{def3.1}
 Let $\rho$ be a distribution on $\R^d$.

(a) Let $d=1$. The distribution $\rho$ is called {\it unimodal} with mode $a$ if\begin{eqnarray}
\rho(dx) = p(x)dx + c\delta_a(dx)
\end{eqnarray}
where $c \ge 0$, $p(x)$ is increasing on $(-\infty,a)$ and decreasing on $(a,\infty)$.

(b) Let $d\geq 1$. Suppose that  $\rho$ is absolutely continuous on $\R^d$. Then $\rho$ is called {\it star-unimodal}  about $0$ if
\begin{eqnarray}
\rho(dx) = p(x)dx 
\end{eqnarray}
where $p(r\xi)$ is decreasing in $r$ on  $(0,\infty)$ for every $\xi \in S^{d-1}.$
\end{defn}

\begin{defn}
\label{def3.2}
 Let $h(x)$ be a  positive measurable function on $\R^d$. 
 The function  $h(x)$ is  called {\it quasi-submultiplicative} on $ \R^d$ if the following conditions are satisfied:

(d.1) For any $R>0$, $\sup_{|x| \le R}h(x) <\infty$

(d.2) For each $\ep >0$, there are $ c_1,  c_2>0$ such that  
\begin{eqnarray}
\label{eq:3.3}
h(x+y)\leq c_1h((1+\ep)x)h(c_2y)\mbox{ for  every } x,  y\in \R^d.
\end{eqnarray}
\end{defn}

\medskip

\begin{rem}
\label{rem3.1}
(i)  Let $h(x)$ be quasi-submultiplicative on $ \R^d$.  
Note that letting $y=0$ in  \eqref{eq:3.3}, we have the following : 

(d.3) For any $\delta>1$, there is $c_3 >0$ such that $h(x) \le c_3
h(\delta x)$ for every $x\in \R^d$.

(ii) Dividing into the two cases $|y| \le \ep |x|$ and $|y| > \ep |x|$ for $\ep >0$,  we see  that if the  function  $h(x)$ is positive and increasing on $\R_+$, then $h(|x|)$  is  quasi-submultiplicative on $ \R^d$. 
\end{rem}

The following celebrated results (i) and (ii) are due to \cite{ya78}  and \cite{wo78}, respectively.

\begin{lem} Suppose that $\mu$ is non-degenerate.
\label{lem3.1}

{\rm (i)} If $d=1$, then $\mu$ is unimodal.

{\rm (ii)}  If $\{X(t)\}$ is symmetric, then $\mu$ is  star-unimodal about $0$.
\end{lem}

The following lemma is from Theorem 25.3 of \cite{sa99}. Refer also to \cite{kr70,kr72} and \cite{sa73}.

\begin{lem}
\label{lem3.2}
  Let $\rho$ be an infinitely divisible distribution on ${\R}^d$ with L\'evy measure $\eta$. Furthermore, let  $h(x)$ be submultiplicative on ${\R}^d$. Then $\int_{{\R}^d}h(x)\rho(dx) <\infty$ if and only if $\int_{|x| >1}h(x)\eta(dx) <\infty.$
\end{lem}

\begin{prop}
\label{prop3.1}
  Let  $h(x)$ be quasi-submultiplicative on ${\R}^d$.

{\rm (i)} If 
\begin{equation}
\label{eq:3.4}
  \int_{{\R}^d}h(x)\mu(dx) <\infty, 
\end{equation}
 then for any $\ep \in (0,1)$ and any $l \in \N$, 
\begin{equation}
\label{eq:3.5}
\int_{{\R}^d}h((1+ \ep)^{-1}x)\rho_l(dx) <\infty. 
\end{equation}

{\rm (ii)} If
\begin{equation}
\label{eq:3.6}
 \int_{{\R}^d}h(cx)\mu(dx) <\infty \mbox{ and } \int_{{\R}^d}h(x)\mu(dx) =\infty
\end{equation}
 for some $c \in (0,\infty)$, then for any $\ep \in (0,1)$, there is $l(\ep)\in\N$ such that for any integer $l\ge l(\ep)$,
\begin{equation}
\label{eq:3.7}
\int_{{\R}^d}h((1+ \ep)x)\rho_l(dx) =\infty.
 \end{equation}
\end{prop}

\prf First we prove (i). Suppose that \eqref{eq:3.4} holds.   Define a subset $\Lambda$ of $(\R^d)^{\N}$ by
\begin{eqnarray*}
 && \Lambda := \{ (x_j)_{j=1}^{\infty} : |x_j | \le N e^{jHl/2} \mbox{ for } j \ge 1 \}.
\end{eqnarray*}
Here we choose $N > 0$ such that $\rho_{l}(\{x : |x| \leq N\}) > 0$. Since the log-moment of $\rho_{l}$ is finite, we have
\begin{eqnarray}
c_0 :=  \int_{\Lambda}\prod_{n=1}^{\infty}\rho_{l}(dx_n) >0.
\end{eqnarray} 
Moreover, since $h(x)$ is quasi-submultiplicative, we see from condition (d.1) of Definition \ref{def3.2}  that
\begin{eqnarray*}
&& c_4:= \sup_{ (x_n)\in \Lambda}h(-(1+\ep)^{-1}c_2\sum_{n=1}^{\infty}e^{-nHl}x_n)< \infty.
\end{eqnarray*} 
Hence we have 
 \begin{eqnarray*}
& &c_0 \int_{\R^d}h((1+\ep)^{-1}x_0)\rho_{l}(dx_0)\\
&& \leq   c_1 \int_{\R^d\times \Lambda}h(\sum_{n=0}^{\infty}e^{-nHl}x_n)h(-(1+\ep)^{-1}c_2\sum_{n=1}^{\infty}e^{-nHl}x_n)\prod_{n=0}^{\infty}\rho_{l}(dx_n)\nonumber \\
&& \leq c_1 c_4 \int_{\R^d\times \Lambda} h(\sum_{n=0}^{\infty}e^{-nHl}x_n) \prod_{n=0}^{\infty}\rho_{l}(dx_n)\nonumber \\
&& \leq c_1 c_4\int_{(\R^d)^{\Z_+}}h(\sum_{n=0}^{\infty}e^{-nHl}x_n)
\prod_{n=0}^{\infty}\rho_{l}(dx_n)\nonumber\\
&& =c_1c_4 \int_{\R^d}h(x)\mu(dx) <\infty.\nonumber
\end{eqnarray*}
Thus we have \eqref{eq:3.5}.

Next we prove (ii).  Let $\ep >0$. Suppose that \eqref{eq:3.6} holds for some $c\in (0,\infty)$. We see that
\begin{eqnarray*}
\qquad E\left(h( X(1))\right)  &\le & c_1 E\left( h((1+\ep)(X(1)-X(e^{-l})))h(c_2 X(e^{-l}))\right)\\
 &= & c_1 E\left(h((1+\ep)(X(1)-X(e^{-l}))\right) E\left(h(c_2 X(e^{-l}))\right)\nonumber \\
 &= & c_1 E\left(h((1+\ep)(X(1)-X(e^{-l}))\right) E\left(h(c_2e^{-Hl} X(1))\right).\nonumber 
\end{eqnarray*}
For all sufficiently large $l$, we have $E\left(h(c_2e^{-Hl} X(1))\right) <\infty$ by (d.3) of Remark \ref{rem3.1}, so we obtain
\begin{eqnarray*}
 E\left(h((1+\ep)(X(1)-X(e^{-l}))\right)= \infty.
\end{eqnarray*}
Hence \eqref{eq:3.7} holds.  \qed
\medskip

Lemmas \ref{lem3.3} and  \ref{lem3.4} below show that it is possible that an infinitely divisible distribution and its L\'evy measure do not have the same tail behaviors. The following is due to \cite{sa73}. Let $B_r$ be the closed ball with center $0$ and radius $r$. Denote by $S(\eta)$ the support of a measure $\eta$ on $\R^d$.

\begin{lem}
\label{lem3.3}
 Let $\rho$ be an infinitely divisible distribution on ${\R}^d$ with L\'evy measure $\eta$. Let $C :=\inf\{r\geq 0: S(\eta)\subset B_r\}$. Then we have the following :

{\rm (i)} If $0 < a <1/C$, then
$ \rho(|x| >r) = o(r^{-ar})$.

{\rm (ii)} If $ a >1/C$, then $ r^{-ar} = o( \rho(|x| >r))$.

{\rm (iii)} We have
\begin{eqnarray}
 \lim_{r \to \infty}\frac{-\log \rho(|x| > r)}{r \log r} = C^{-1}.
\end{eqnarray}  
\end{lem}

\medskip

Let $q(t)$ be a positive and right-continuous function on $(0,\infty)$. Define the function $\varphi(r)$ by $\log r=\int_0^{\varphi(r)}q(t)dt$. If $q(t)=1$ on $\R_+$, then $\varphi(r)=\log r$. The limit of $-\log\rho(|x|>r))/(r\varphi(r))$ is discussed by \cite{ka82} in the case where $q(t)$ is increasing:

{\begin{lem}
\label{lem3.4}
 Let $\rho$ be an infinitely divisible distribution on $\R^d$ with L\'evy measure $\eta$. Suppose that $h$ and $f$ are Young conjugate functions. Let $0 <\gamma \le\infty$. Then we have
\begin{eqnarray}
 \liminf_{r \to \infty}\frac{-\log \rho(|x| > r)}{r  f^{-1}(\log r)} = \gamma
\end{eqnarray}  
if and only if 
\begin{eqnarray}
 \liminf_{r \to \infty}\frac{h^{-1}(-\log \eta(|x| > r))}{r } = \gamma.
\end{eqnarray}  
\end{lem}

We give a sufficient condition for which an infinitely divisible distribution and its L\'evy measure  have the same tail behaviors in the relation ``$\asymp$''. The result is found also in \cite{sw05} in the case of an infinitely divisible distribution on ${\R}_+$. It is proved in Sect.5. Recall the definition \eqref{eq:2.1} of $\widetilde{\eta}$ for a measure $\eta$.

\begin{thm}
\label{th3.1}
  Let $\rho$ be an infinitely divisible distribution on ${\R}^d$ with L\'evy measure $\eta$.
Define a distribution $\zeta$ on $\R_+$ by  $\zeta(x >r):= \widetilde\eta(|x| >r)$ for $r \ge 0.$

{\rm (i)} $\rho\in {\mathcal O}{\mathcal S}$  on ${\R}^d$ if and only if there is a positive integer $n$ such that $\widetilde\eta^{n*}\in {\mathcal O}{\mathcal S}$ on $\R^d$. Moreover, if $\rho\in {\mathcal O}{\mathcal S}$, then, for the above $n$,
\begin{eqnarray}
\label{eq:3.12}
\rho(|x| > r) \asymp \widetilde\eta^{n*}(|x| >r).
\end{eqnarray}

{\rm (ii)} If $\widetilde\eta  \in {\mathcal O}{\mathcal S}$, then we have
\begin{eqnarray}
\label{eq:3.13}
\rho(|x| > r) \asymp \eta(|x| >r).
\end{eqnarray}

{\rm (iii)}  $\rho\in {\mathcal O}{\mathcal S}$ if and only if there is a positive integer $m$ such that $\zeta^{m*}\in {\mathcal O}{\mathcal S}$ on $\R_+$ and 
\begin{eqnarray}
\label{eq:3.14}
\rho(|x| > r) \asymp \zeta^{m*}(x >r).
\end{eqnarray}

{\rm (iv)} Let $h(x)$ be a nonnegative increasing function on $\R_+$. Suppose that $\widetilde\eta \in {\mathcal O}{\mathcal S}$ on $\R^d$. Then  $\int_{{\R}^d}h(|x|)\rho(dx) <\infty$ if and only if $\int_{|x| >1}h(|x|)\eta(dx) <\infty.$
\end{thm}

\begin{rem}
\label{rem3.2}
 In (i) and (iii) of Theorem \ref{th3.1}, $n$ and $m$ are not always 1. See Theorem 1.1 of \cite{sw05}. In (iii), we do not know whether  $\rho\in {\mathcal O}{\mathcal S}$ provided that there is a positive integer $m \ge 2$ such that $\zeta^{m*}\in {\mathcal O}{\mathcal S}$ on $\R_+$.
\end{rem}

We give a result as to the class ${\mcal D}$. The class ${\mcal OS}$ includes the class ${\mcal D}$. The fact is found also in \cite{wa96} in the case of infinitely divisible distributions on $\R_+$. 

\begin{cor}
\label{cor3.1}
 Let $\rho$ be an infinitely divisible distribution on ${\R}^d$ with L\'evy measure $\eta$. Then the following holds:

{\rm (i)} $ \rho \in {\mcal D}$ if and only if   $\widetilde\eta \in {\mcal D}$.  

{\rm (ii)} If  $ \widetilde\eta  \in {\mcal D}$, then $\rho(|x|>r)\asymp \eta(|x|>r)$.
\end{cor}

\prf Let $\zeta$ be the distribution on $\R_+$ defined in Theorem \ref{th3.1}. Note that ${\mcal D} \subset {\mathcal O}{\mathcal S}$ on $\R_+$ and on $\R^d$. Thus assertion (ii) is obvious from  Theorem \ref{th3.1} (ii). 
Next we prove assertion (i). See Lemma \ref{lem5.3} (ii) in Sect.5. We find from Theorem \ref{th3.1} (i) that  if $\widetilde\eta \in {\mcal D}$, then $\rho \in {\mcal D}$. 
We obtain from  Theorem \ref{th3.1} (iii) that if $\rho \in {\mcal D}$, then   $\zeta^{m*}\in {\mathcal D}$ on $\R_+$ for some $m \in \N$.  Moreover we see from Proposition 1.1 (iii) and 2.5 (iii)  of \cite{sw05} (Lemma \ref{lem5.4} in Sect.5) that $\zeta^{m*}\in {\mathcal D}$ on $\R_+$ for some $m\in \N$ if and only if $\zeta\in {\mathcal D}$ on $\R_+$.  Thus, if $\rho \in {\mcal D}$, then   $\zeta\in {\mathcal D}$ and thereby assertion (i) is true.   \qed

\begin{prop}
\label{prop3.2}
{\rm (i)} $G_1(r)\in  OR$ if and only if $L(r) \in  OR$.

{\rm (ii)} If $L(r) \in  OR$, then $ G_1(r) \asymp L(r).$
\end{prop}

\prf Let $b:= e^{-H}$ and $\delta :=e^H-1$. On the one hand, we see 
\begin{eqnarray*}
\mu(|x| > r) &=& P(|X(e^{-1})+X(1)-X(e^{-1})| > r)\\
&\le& P(e^{-H}|X(1)| > \sqrt{b} r) + P(|X(1)-X(e^{-1})| > (1-\sqrt{b})r)\\
&=& \mu(|x|>r/\sqrt{b})+\rho_1(|x|>(1-\sqrt{b})r).
\end{eqnarray*}
Hence we have
\begin{eqnarray*}
\mu(r < |x| \le r/\sqrt{b}) \le  \rho_1(|x| >(1-\sqrt{b})r). 
\end{eqnarray*}
We obtain that
\begin{eqnarray}
\label{eq:3.15}
\mu(r < |x| \le e^{3H}r) &\le&\sum_{j=0}^5\mu(b^{-j/2}r < |x| \le b^{-(j+1)/2}r)\nonumber\\
 &\le &\sum_{j=0}^5\rho_1(|x| > (1-\sqrt{b})b^{-j/2}r)\nonumber \\
&\le& 6 \rho_1(|x| >(1-\sqrt{b})r).  
\end{eqnarray}
On the other hand, we see
\begin{eqnarray*}
& &\mu(|x| > r)  = P(|X(e^{-1})+X(1)-X(e^{-1})| > r)\\
&&\ge P(|X(1)-X(e^{-1})| >e^H r,|X(e^{-1})| \le \delta r)\\
&& \quad +P(|X(e^{-1})+X(1)-X(e^{-1})| >r, |X(1)-X(e^{-1})| \le e^H r)\\
&&\ge P(|X(1)-X(e^{-1})| >e^H r)P(|X(e^{-1})| \le \delta r)\\
&&\quad +P(|X(e^{-1})| >e^{2H} r)P(|X(1)-X(e^{-1})| \le e^H r)\\
&&=\rho_1(|x|>e^H r)\mu(|x|\leq \delta r/b)+\mu(|x|>e^{3H}r)\rho_1(|x|\leq e^H r).
\end{eqnarray*}
Hence we have
\begin{eqnarray}
\label{eq:3.16}
&&\mu(r <|x| \le e^{3H} r)\ge \rho_1(|x|>e^H r)\left\{\mu(|x| \le \delta r/b)- \mu(|x| >e^{3H}r)\right\}. \qquad 
\end{eqnarray}
By \eqref{eq:3.15} and \eqref{eq:3.16}, $L(r): =\mu(r < |x| \le e^{3H} r) \in   OR$
if and only if $G_1(r): =\rho_1(|x| >r) \in  OR$. Furthermore, if $L(r) \in  OR$, then  $L(r) \asymp G_1(r). $ \qed

\begin{thm}
\label{th3.2}

{\rm (i)} The following are equivalent:

{\rm (1)}   $G_1(r)\in OR,$\quad {\rm (2)} $K(r) \in OR$,\quad {\rm (3)} $L(r) \in  OR$.

{\rm (ii)} If $K(r) \in  OR$, then $G_1(r) \asymp K(r) \asymp L(r)$.
\end{thm}

\prf Now \eqref{eq:2.11} implies that $\eta_1(|x|>r)\in OR$ is equivalent to $K(r)\in OR$. Use Corollary \ref{cor3.1} with $\rho=\rho_1$ and $\eta=\eta_1$. Hence $\eta_1(|x|>r)\in OR$, namely, $\widetilde{\eta}_1\in{\mathcal D}$ is equivalent to $\rho_1\in{\mathcal D}$, namely, $G_1(r):=\rho_1(|x|>r)\in OR$. Furthermore, if $K(r)\in OR$, namely, $\widetilde{\eta}_1\in{\mathcal D}$, then $G_1(r)=\rho_1(|x|>r)\asymp \eta_1(|x|>r) \asymp K(r)$. The remaining proof is proved by Proposition \ref{prop3.2}.\qed

\begin{prop} 
\label{prop3.3}
Suppose either that $d=1$ or that $\{X(t)\}$ is symmetric and non degenerated with $d \ge 2$. 

{\rm (i)} $M(r) \in  OR$ if and only if $L(r) \in  OR$.

{\rm (ii)} If $M(r) \in  OR$, then $r^dM(r) \asymp L(r)$.
\end{prop}

\prf We have 
$$L(r):=\mu(r<|x|\leq e^{3H}r)=\int_r^{e^{3H}r}M(u)u^{d-1}du,$$
By Lemma \ref{lem3.1}, the proposition is obvious. \qed 



\section{Proof of the results}

In this section, we prove the results mentioned in Sect.2. First of all, we present two important propositions which lead to the key theorem, that is, Theorem \ref{th2.1}.

\begin{prop}
\label{prop4.1}
Let $g\in{\mathcal G}_1$ and $l\in \N$.

{\rm (i)}\ If
\begin{eqnarray}
\int_0^\infty P(|X(1)-X(e^{-l})|>g(x))dx = \infty,
\end{eqnarray}
then 
\begin{eqnarray}
\label{eq:4.2}
\limsup_{t\to\infty}\frac{|X(t)|}{t^Hg(\log t)} \geq 1\qquad \mbox{a.s.}
\end{eqnarray}

{\rm (ii)}\ If
\begin{eqnarray}
\label{eq:4.3}
\int_0^\infty P(|X(1)-X(e^{-l})|>g(x))dx < \infty,
\end{eqnarray}
then 
\begin{eqnarray}
\label{eq:4.4}
\limsup_{t\to\infty}\frac{|X(t)|}{t^Hg(\log t)} \leq (1-e^{-lH})^{-1}\qquad \mbox{a.s.}
\end{eqnarray}
\end{prop}

\medskip

\begin{prop}
\label{prop4.2}
Let $g\in{\mathcal G}_1$ and $l\in \N$. Suppose that ${\displaystyle \limsup_{x\to\infty}g(x+1)/g(x)<e^H}$.

{\rm (i)}\ If
\begin{eqnarray}
\int_0^\infty P(|X(1)-X(e^{-l})|>g(x))dx = \infty,
\end{eqnarray}
then 
\begin{eqnarray}
\limsup_{t\to 0}\frac{|X(t)|}{t^Hg(|\log t|)} \geq 1\qquad \mbox{a.s.}
\end{eqnarray}

{\rm (ii)}\ If
\begin{eqnarray}
\label{eq:4.7}
\int_0^\infty P(|X(1)-X(e^{-l})|>g(x))dx < \infty,
\end{eqnarray}
then 
\begin{eqnarray}
\label{eq:4.8}
\limsup_{t\to 0}\frac{|X(t)|}{t^Hg(|\log t|)} \leq (1-e^{-lH_1})^{-1}\qquad \mbox{a.s.}
\end{eqnarray}
for any $H_1\in(0,H)$. 
\end{prop}

\medskip

In order to prove Propositions \ref{prop4.1} and \ref{prop4.2}, we use some results of \cite{wa02} for shift selfsimilar additive random sequences defined below.

\begin{defn}
\label{def4.1}
  Let $c > 1$.  An $\R^d$-valued random sequence $\{Y(n), n \in \Z \}$ on a probability space $({\widetilde \Om},
{\widetilde {\mathcal{F}}},{\widetilde P})$ is called a  {\em shift $c$-self-similar additive random sequence} if 
the following two conditions are satisfied:

(1)  The sequence $\{Y(n),  n \in \Z \}$ has shift $c$-self-similarity, that is,
$$ \{Y(n+1),  n \in \Z  \}\stackrel{\rm d}{ =} \{cY(n),  n \in \Z \},  $$
where the symbol $\stackrel{\rm d}{ =}$ stands for equality in   the finite-dimensional distributions.

(2) The sequence $\{Y(n),  n \in \Z \}$ has independent increments (or additivity), that is, for every $ n \in \Z$, $\{Y(k), k \leq n\}$ and  $Y(n+1)-Y(n)$ are independent.
\end{defn}

In the rest of this section, we define the random sequence $\{Y(n), n \in \Z \}$ by $Y(n):= X(e^n)$ for $n \in \Z$. Under the assumption that $\{X(t)\}$ is not purely Gaussian, we decompose  $\{X(t)\}$ in law as the sum of two independent selfsimilar additive processes  $\{X_j(t)\}$ for $j =1,2$ as follows: Denote by $\mu_j$ the distribution of $X_j(1)$ for $j =1,2$. Define $k_{(\xi, N)} (r) : =k_{\xi} (r\vee N^{-1})$ for  $N >0$ and denote by $ K_1(r) $ the function $K(r)$ in \eqref{eq:2.8} replacing $k_{\xi}(r)$ by $k_{(\xi, N)} (r)$. Here we note that  $K (r)- K_1(r)=0$ for $r \ge N^{-1}$. By Lemma \ref{lem2.2}, we can define $\{X_j(t)\}$ by determining $\mu_j$ for $j =1,2$. Now $\mu$ satisfies \eqref{eq:2.2}, \eqref{eq:2.3} and \eqref{eq:2.5}. Hence we take $\mu_1$ and $\mu_2$ such that $\mu = \mu_1* \mu_2$ and for $j=1,2$,
 
\begin{equation}
 \widehat \mu_j(z):= \int_{\R^d}\exp (i\la z, x\ra)\mu_j(dx) = \exp(\ps_j(z)),\qq z\in\R^d
\end{equation}
with
\begin{equation}
\ps_1(z):=\int_{S^{d-1}}\sg(d\xi)\int_0^{\infty}(e^{i\la z,r\xi\ra}-1) k_{(\xi,N)}(r)r^{-1}dr,
\end{equation}
and
\begin{equation}
\ps_2(z):=\ps(z) -\ps_1(z).
\end{equation}

\medskip

The following is obvious from the definition.

\medskip

\begin{lem}
\label{lem4.1}
 The sequence  $\{Y(n), n \in \Z \}$ 
is a shift $e^H$-self-similar additive random sequence.
\end{lem}

The following lemma is from Lemma 5.2 in \cite{wa02}. Remark \ref{rem4.1} is
needed in the proof of Lemma \ref{lem4.2} (i) below. 

\begin{rem}
\label{rem4.1}
 Let $g \in {\cal G}_1 $ and let $l \in \N$. If
\begin{equation}
 \int_0^{\infty}P(|Y(0)-Y(-l)| >g(x) )dx = \infty, 
\end{equation}
then, for all $k \in {\bf N}$ and all $\ep \in (0,1)$,  
\begin{equation}
 \int_0^{\infty}P(|Y(0)-Y(-kl)| >(1-\ep )g(x) )dx = \infty. 
\end{equation}
\end{rem}

\medskip

\begin{lem}
\label{lem4.2}
  Let $g \in {\cal G}_1 $ and let $l \in \N$. 

{\em (i) } If
\begin{equation}
 \int_0^{\infty}P(|Y(0)-Y(-l)| >g(x) )dx = \infty,\label{lem4.2.1}
\end{equation}
then  
\begin{equation}
 \limsup_{n \to \infty}\frac{|Y(n)|}{e^{nH}g(n)} \geq 1 \quad\mbox{ {\em a.s.}}
\end{equation} 

{\em (ii) } Suppose that ${\displaystyle \limsup_{x\to\infty}g(x+1)/g(x)<e^H}$. If \eqref{lem4.2.1} is satisfied, then  
\begin{equation}
 \limsup_{n \to \infty}\frac{|Y(-n)|}{e^{-nH}g(n)} \geq 1 \quad \mbox{ {\em a.s.}}  
\end{equation}
\end{lem}

\medskip

In Lemmas \ref{lem4.3} and \ref{lem4.4} below, we assume that $\{X(t)\}$ is non-Gaussian. Further, we take sufficiently large $N$ and let $\{X_1(t)\}$ be non-Gaussian.

\medskip

\begin{lem}
\label{lem4.3}
 Let $l \in \N.$  There is $c_1>0$ such that, for $r \ge 0$,
\begin{equation}
\sup_{t \in [e^{-l},1]} P(|X_1(t)-X_1(e^{-l})| >r) \le c_1 P(|X_1(1)-X_1(e^{-l})| >r).\label{lem4.3.1}
\end{equation}
\end{lem}

\prf Note that the distribution of $X_1(t) -X_1(e^{-l})$ for $t \in (e^{-l},1]$ is compound Poisson. Denote by $\nu_t$ its L\'evy measure. Then $\nu_t$ is represented as
\begin{eqnarray*}
\nu_t(B)=\int_{S^{d-1}}\sg(d\xi)\int_0^\infty 1_{B}(r\xi)\frac{k_{(\xi,N)}(r/t^H)-k_{(\xi,N)}(e^{lH}r)}{r}dr
\end{eqnarray*}
for any Borel set $B$ in $\R^d$. We have, for $r \ge 0$,
\begin{eqnarray*}
P(|X_1(t)-X_1(e^{-l})| >r) = e^{-c(t)} \sum_{n=1}^{\infty}\frac{1}{n!}\nu_t^{n*}(|x| >r),
\end{eqnarray*}
where $c(t):= \nu_t(\R^d)< \infty.$
Then $\nu_t(B)$ is increasing in $t$, so is $\sum_{n=1}^{\infty}\nu_t^{n*}(B)/n!$. Thus we obtain that 
\begin{equation*}
\sup_{t \in [e^{-l},1]} P(|X_1(t)-X_1(e^{-l})| >r) \le e^{c(1)}P(|X_1(1)-X_1(e^{-l})| >r).
\end{equation*}
Thus we have \eqref{lem4.3.1}. \qed

\begin{lem} 
\label{lem4.4}
Let $\delta >0$ and $l\in \N$. If  $N$ of $k_{(\xi,N)}(r)$ is sufficiently large, we can take sufficiently large $\delta_0 >\delta $ such that
\begin{equation}
\sup_{t \in [e^{-l},1]} P(|X_2(t)-X_2(e^{-l})| >r) = o(\exp(-\delta_0 r\log r)) \quad \mbox{as $r \to \infty$}.\label{lem4.4.1}
\end{equation}
\end{lem}

\prf We see that
\begin{eqnarray*}
&&\sup_{t \in [e^{-l},1]} P(|X_2(t)-X_2(e^{-l})| >r)\\
&&\le \sup_{t \in [e^{-l},1]}( P(|X_2(t)| >r/2)+ P(|X_2(e^{-l})| >r/2))\\
&& \le 2 P(|X_2(1)| >r/2).
 \end{eqnarray*}
In the last inequality, we used the selfsimilarity.
Now the support of the L\'evy measure of $\{X_2(t)\}$ is the empty set or is contained in a small ball with center 0. Hence we obtain \eqref{lem4.4.1} from Lemma \ref{lem3.3}. \qed

\begin{lem}
\label{lem4.5}
 Let $l\in\N$. For any $\ep \in (0,1)$, there is  $c_2>0$ such that, for $r \ge 0$,
\begin{equation}
\sup_{t \in [e^{-l},1]} P(|X(t)-X(e^{-l})| >r) \le c_2 P(|X(1)-X(e^{-l})| >(1-\ep)r).\label{lem4.5.1}
\end{equation}
\end{lem}

\medskip

\prf  In the case where $\{X(t)\}$ is Gaussian, the lemma can be
proved by straightforward calculations. Thus we only prove the case
where $\{X(t)\}$ is non-Gaussian. Choose $\delta \in (0,1)$ such that
$1-\ep= (1-\delta)^2$. By Lemmas \ref{lem4.3} and \ref{lem4.4}, there is $c>0$ such that
for sufficiently large $N$ and $\delta_0$,   
\begin{eqnarray}
&&\sup_{t \in [e^{-l},1]} P(|X(t)-X(e^{-l})| >r)\nonumber\\
&\le &\sup_{t \in [e^{-l},1]} P(|X_1(t)-X_1(e^{-l})|+ |X_2(t)-X_2(e^{-l})|>r)\nonumber\\
&\le &\sup_{t \in [e^{-l},1]} P(|X_1(t)-X_1(e^{-l})| >(1-\delta)r)+\!\!\!\sup_{t \in [e^{-l},1]}P (|X_2(t)-X_2(e^{-l})|>\delta r)\nonumber\\
&\le & c P(|X_1(1)-X_1(e^{-l})| >(1-\delta)r)+ o(\exp(-\delta_0 \delta r\log (\delta r))).\label{lem4.5.2}
\end{eqnarray}
Further we have
\begin{eqnarray}
 && P(|X_1(1)-X_1(e^{-l})| >(1-\delta)r)\nonumber\\
&\le& P(|X(1)-X(e^{-l})|+ |X_2(1)-X_2(e^{-l})|>(1-\delta)r)\nonumber\\
&\le& P(|X(1)-X(e^{-l})|>(1-\delta)^2 r)+P( |X_2(1)-X_2(e^{-l})|>\delta(1-\delta)r)\nonumber\\
&=&  P(|X(1)-X(e^{-l})|>(1-\ep) r)+ o(\exp(-\delta_0 \delta(1-\delta)
r\log (\delta (1-\delta)r))).\nonumber\\
&&\label{lem4.5.3}
\end{eqnarray}
Taking sufficiently large $\delta_0$, we find from Lemma \ref{lem3.3} that 
\begin{eqnarray*}
 \exp(-\delta_0 \delta r\log (\delta r))) = o(P(|X(1)-X(e^{-l})|>(1-\ep) r))
\end{eqnarray*}
and
\begin{eqnarray*}
\exp(-\delta_0 \delta(1-\delta) r\log (\delta (1-\delta)r)))= o(P(|X(1)-X(e^{-l})|>(1-\ep) r)).\nonumber
\end{eqnarray*}
Thus \eqref{lem4.5.2} and \eqref{lem4.5.3} yields that \eqref{lem4.5.1} holds for some $c_2>0$. \qed

\begin{lem}
\label{lem4.6}
 Let $0 \le s < t$, $a , b >0$, and $\ep >0$. Then
\begin{equation}
 P(\sup_{u \in [s,t]}|X(u)-X(s)| >3\ep) \le 3 \sup_{u \in [s,t]}P(|X(u)-X(s)| >\ep),\label{lem4.6.1}
\end{equation}
and
\begin{equation}
 P(\sup_{u \in [s,t]}|X(u)-X(s)| >a+b) \le \frac{P(|X(t)-X(s)| >a)}{P( \sup_{u \in [s,t]}|X(u)-X(s)| \le b/2)}.\label{lem4.6.2}
\end{equation}
\end{lem}

\prf The inequality \eqref{lem4.6.1} is due to Lemma 2.1 of \cite{ya03}. The proof of \eqref{lem4.6.2} follows along the lines of the proof of \eqref{lem4.6.1} from Remark 20.3 in \cite{sa99}.
\qed

\begin{lem}
\label{lem4.7}
 For any $\ep \in (0,1)$, there is    $R >0$ such that, for $r > R$,
\begin{equation}
 P(\sup_{t \in [e^{-l},1]}|X(t)-X(e^{-l})| >r) \le 2P(|X(1)-X(e^{-l})| >(1-\ep)r).\label{lem4.7.1}
\end{equation}
\end{lem}

\medskip

\prf We see from \eqref{lem4.6.2} of Lemma \ref{lem4.6} that
\begin{eqnarray*}
P(\sup_{t \in [e^{-l},1]}|X(t)-X(e^{-l})| >r) &\le& \frac{P(|X(1)-X(e^{-l})| >(1-\ep)r)}{P( \sup_{t \in [e^{-l},1]}|X(t)-X(e^{-l})| \le \ep r/2)}.\nonumber
\end{eqnarray*}
We obtain from Lemma \ref{lem4.5} and \eqref{lem4.6.1} of Lemma
\ref{lem4.6} that
\begin{eqnarray*}
P( \sup_{t \in [e^{-l},1]}|X(t)-X(e^{-l})| \le \ep r/2) &=& 1-P( \sup_{t \in [e^{-l},1]}|X(t)-X(e^{-l})| > \ep r/2)\nonumber\\
 &\ge& 1-3\sup_{t \in [e^{-l},1]}P( |X(t)-X(e^{-l})| > \ep r/6)\nonumber\\
 &\ge&1-3c_2P( |X(1)-X(e^{-l})| >(1-\ep) \ep r/6).\nonumber
\end{eqnarray*}
Taking sufficiently large $R >0$ such that 
$$3c_2P( |X(1)-X(e^{-l})| >(1-\ep) \ep r/6) < 1/2 \q \mbox{ for } r >R,$$
 we have  \eqref{lem4.7.1}. \qed

\medskip

\noindent
{\it Proof of Proposition \ref{prop4.1}:} Assertion (i) follows from (i) of Lemma \ref{lem4.2}. Next we prove (ii). Suppose that \eqref{eq:4.3} holds.
Let $M>0$ be a sufficiently large positive integer. Take a positive constant $b$ such that $e^l\geq b >1$ and $b^M=e$. We see from Lemma \ref{lem4.7} that, for any $\ep \in (0,1)$, there is $R >0$ such that, for $r > R$,
\begin{equation*}
 P(\sup_{t \in [b^{-1},1]}|X(t)-X(e^{-l})| >r) \le 2P(|X(1)-X(e^{-l})| >(1-\ep)r).
\end{equation*}
Thus, by using the selfsimilarity, we have
\begin{eqnarray*}
& & P(\sup_{t \in [b^{-m-1}e^{n},b^{-m}e^n]}|X(t)-X(b^{-m}e^{n-l})| >(b^{-m}e^n)^Hr)\\
& \le&  2P(|X(1)-X(e^{-l})| >(1-\ep)r)\nonumber
\end{eqnarray*}
for $0 \le m \le M-1$ and $n \in \Z$. Let $r_n:= g(n-1)/(1-\ep)$. We obtain from \eqref{eq:4.3} that
\begin{eqnarray*}
\sum_{n=1}^{\infty} P(\sup_{t \in [b^{-m-1}e^{n},b^{-m}e^n]}|X(t)-X(b^{-m}e^{n-l})| >(b^{-m}e^n)^Hr_n) <  \infty.
\end{eqnarray*}
By virtue of the Borel-Cantelli lemma, there is a random number $N_m$ such that, for any $n \ge N_m$, almost surely
\begin{eqnarray}
\label{eq:4.25}
&& \sup_{t \in  [b^{-m-1}e^{n},b^{-m}e^n]}|X(t)|-|X(b^{-m}e^{n-l})| \le \frac{(b^{-m}e^n)^Hg(n-1)}{1-\ep}. \quad
\end{eqnarray}
Let $k:=k_{(m,n)}$ be the maximum integer satisfying $ n - k_{(m,n)}l \ge N_m$. Let  $1 \le j \le k$. Then we substitute $n-jl$ for $n$ in \eqref{eq:4.25} and   see that, for $1 \le j \le k$, almost surely
\begin{eqnarray}
\label{eq:4.26}
&& |X(b^{-m}e^{n-jl})|-|X(b^{-m}e^{n-(j+1)l})| \le \frac{(b^{-m}e^{n-jl})^Hg(n-jl-1)}{1-\ep}.\qquad
\end{eqnarray}
Recall that $e=b^M$. Adding up \eqref{eq:4.25} and \eqref{eq:4.26} for $1 \le j \le k$, we have
\begin{eqnarray*}
\sup_{t \in [b^{Mn-m-1},b^{Mn-m}]}|X(t)|-|X(b^{-m}e^{n-(k+1)l})|& \le& \sum_{j=0}^k\frac{(b^{-m}e^{n-jl})^Hg(n-jl-1)}{1-\ep}\nonumber\\
& \le& \frac{(b^{Mn-m})^Hg(n-1)}{(1-\ep)(1-e^{-lH})}.\nonumber
\end{eqnarray*}
Notice that $N_m-l\leq n-(k+1)l<N_m$. Thus,
\begin{eqnarray*}
\limsup_{t\to\infty}\frac{|X(t)|}{t^Hg(\log t)} & \le &\limsup_{n\to\infty}\sup_{t \in [b^{M(n-1)},b^{Mn}]}\frac{ |X(t)|}{t^Hg(\log t)} \nonumber\\
& \le &\limsup_{n\to\infty}\sup_{0 \le m\le M-1}\sup_{t \in [b^{Mn-m-1},b^{Mn-m}]}\frac{b^H |X(t)|}{(b^{Mn-m})^Hg(n-1)}\nonumber\\
& \le & \!\!\frac{b^H}{(1-\ep)(1-e^{-lH})}+\limsup_{n\to\infty}\!\max_{N_m-l\leq i<N_m}\frac{b^H|X(b^{-m}e^i)|}{(b^{Mn-m})^Hg(n-1)}\\
&=& \frac{b^H}{(1-\ep)(1-e^{-lH})}.
\end{eqnarray*}
Letting $\ep\downarrow 0$ and $b\downarrow 1$, we have \eqref{eq:4.4}. \qed
\medskip

\noindent
{\it Proof of Proposition \ref{prop4.2}:} Assertion (i) follows from (ii) of Lemma \ref{lem4.2}. Next we prove (ii). Suppose that \eqref{eq:4.7} holds. In the same way as in the proof of Proposition \ref{prop4.1}, there is a random number $N$ such that, for any $n\geq N$, almost surely
\begin{eqnarray*}
&& \sup_{t \in  [e^{-l(n+1)},e^{-ln}]}|X(t)-X(e^{-(n+1)l})| \le \frac{e^{-lnH}g((n+1)l)}{1-\ep}. 
\end{eqnarray*}
Let $0<H'<H$. Then we have $g(n+1)/g(n)<e^{H'}$ for all sufficiently large $n$. Put $H_1=H-H'$. Note from Theorem 2.1 of \cite{wa02} that ${\displaystyle \lim_{n\to \infty}X(e^{-ln})=0}$ a.s. Hence we have 
\begin{eqnarray*}
\sup_{t \in  [e^{-l(n+1)},e^{-ln}]}|X(t)| &\le& (1-\epsilon)^{-1}\sum_{j=0}^\infty e^{-l(n+j)H}g(l(n+j+1))\\
&\le& (1-\epsilon)^{-1}e^{-lnH}\sum_{j=0}^\infty e^{-ljH_1}g(l(n+1))\\
&=& \frac{e^{-lnH}g(l(n+1))}{1-\ep}\cdot(1-e^{-lH_1})^{-1}
\end{eqnarray*}
for all sufficiently large $n$. Here we have 
\begin{eqnarray*}
\frac{e^{-lnH}g(l(n+1))}{t^Hg(|\log t|)}\leq \frac{g(l(n+1))}{e^{-lH}g(nl)}<1 
\end{eqnarray*}
for $t\in[e^{-l(n+1)},e^{-ln}]$. Therefore,
\begin{eqnarray*}
\limsup_{t\to 0}\frac{|X(t)|}{t^Hg(|\log t|)}\leq (1-\ep)^{-1}\cdot(1-e^{-lH_1})^{-1}.
\end{eqnarray*}
Letting $\ep\to 0$, we obtain \eqref{eq:4.8}.\qed
\medskip

\noindent
{\it Proof of Theorem \ref{th2.1}:} Note that, for $\delta >0$ and $l \in \N$,
\begin{eqnarray}
\label{eq:4.27}
\int_0^\infty P(|X(1)-X(e^{-l})|>\delta g(x))dx &=& \int_0^{\infty} G_l(\delta g(x))dx\nonumber\\
& =&\int_{\R^d} g^{-1}(\delta^{-1}|x|)\rho_l(dx).
\end{eqnarray}
Thus the proof of  (ii) is clear from Proposition \ref{prop4.1}.  Assertion (i) follows from (ii).
 Suppose that $g \in{\mathcal G}_1$ and $g^{-1}(|x|)$ is submultiplicative on $\R^d$. We prove the equivalence of \eqref{eq:2.13} and \eqref{eq:2.14}. Recall the L\'evy measure $\eta_1$ of $\rho_1$, that is, \eqref{eq:2.11}.
We find from Lemma \ref{lem3.2}  that
\begin{eqnarray}
\label{eq:4.28}
\int_1^\infty \frac{K(r)-K(e^Hr)}{r}g^{-1}(\delta^{-1} r)dr < \infty
\end{eqnarray}
if and only if  $\int_{{\R}^d}g^{-1}(\delta^{-1} |x|)\rho_1(dx) <\infty.$
For some $c_1>0$, we have
\begin{eqnarray*}
\int_{\R^d}g^{-1}(\delta^{-1}|x|)\rho_l(dx)& =& \int_{(\R^d)^l}g^{-1}\left(\delta^{-1}\left|\sum_{k=1}^le^{(-k+1)H}x_k\right|\right)\prod_{k=1}^l\rho_1(dx_k). \\
& \le &c_1^l\int_{(\R^d)^l}\prod_{k=1}^lg^{-1}(\delta^{-1}e^{(-k+1)H}|x_k|)\rho_1(dx_k)\nonumber\\
& \le &c_1^l\left(\int_{\R^d}g^{-1}(\delta^{-1}|x|)\rho_1(dx)\right)^l. \nonumber
\end{eqnarray*}
Thus \eqref{eq:4.28} is equivalent to that $\int_{\R^d}g^{-1}(\delta^{-1}|x|)\rho_l(dx) <\infty$ for all $l$. Suppose that \eqref{eq:2.14} holds. 
Let $\delta > C$. Then  \eqref{eq:4.28} holds, and thus $\int_{\R^d}g^{-1}(\delta^{-1}|x|)\rho_l(dx) <\infty$ for all $l$. Let  $\delta < C$. Then we have
\begin{eqnarray*}
\infty=\int_1^\infty \frac{K(r)-K(e^Hr)}{r}g^{-1}(\delta^{-1} r)dr\le\int_1^\infty \frac{K(r)-K(e^{lH}r)}{r}g^{-1}(\delta^{-1} r)dr.\nonumber
 \end{eqnarray*}
By Lemma \ref{lem3.2},  $\int_{\R^d}g^{-1}(\delta^{-1}|x|)\rho_l(dx) =\infty$ for all $l$, and thereby \eqref{eq:2.13} holds. Conversely, suppose that \eqref{eq:2.13} holds. Let $\delta > C$. 
 Then  $\int_{{\R}^d}g^{-1}(\delta^{-1} |x|)\rho_l(dx) <\infty$ for all $l$, and hence \eqref{eq:4.28} holds. Let  $\delta < C$. Since $\int_{\R^d}g^{-1}(\delta^{-1}|x|)\rho_l(dx) =\infty$ for some $l$, we see that
\begin{eqnarray*}
\infty& = &\int_1^\infty \frac{K(r)-K(e^{lH}r)}{r}g^{-1}(\delta^{-1} r)dr\\ 
&= &\sum_{k=1}^{l}\int_1^\infty \frac{K(e^{(k-1)H}r)-K(e^{kH}r)}{r}g^{-1}(\delta^{-1} r)dr\nonumber\\
&= &\sum_{k=1}^{l}\int_{e^{(k-1)H}}^\infty \frac{K(r)-K(e^Hr)}{r}g^{-1}(\delta^{-1}e^{-(k-1)H} r)dr\nonumber\\
& \le & l \int_1^\infty \frac{K(r)-K(e^Hr)}{r}g^{-1}(\delta^{-1} r)dr.  \nonumber
 \end{eqnarray*}
Hence \eqref{eq:2.14} holds. We have   proved  the equivalence of
\eqref{eq:2.13} and
\eqref{eq:2.14}. Thus assertion (iii) follows from (ii). The second assertion
of (iii) is trivial. \qed
\medskip

\noindent
{\it Proof of Theorem \ref{th2.2}:} We see from Theorem \ref{th3.2} that $F(r) \in OR$ is equivalent to  $G_1(r) \in OR$.  Suppose that $G_1(r) \in OR$ and there is $g(x) \in {\cal G}_1 $ such that \eqref{eq:1.2} holds with $C=1$. Then we see from Proposition \ref{prop4.1}  that
\begin{eqnarray*}
&& \int_0^{\infty} G_1(2^{-1}(1-e^{-H})g(x) )dx\\
&&=\int_0^{\infty}P(|X(1)-X(e^{-1})| >2^{-1}(1-e^{-H})g(x) )dx = \infty 
\end{eqnarray*}
and
\begin{eqnarray*}
&& \int_0^{\infty} G_1(2g(x) )dx=\int_0^{\infty}P(|X(1)-X(e^{-1})| >2g(x) )dx < \infty. 
\end{eqnarray*}
As $G_1(r)\in OR$, this is a contradiction. Hence if $G_1(r)\in OR$, then there is no $g(x) \in {\cal G}_1 $ such that \eqref{eq:1.2} holds with $C=1$. Conversely, suppose that $G_1(r)  \notin OR$.
 Then there is a positive sequence $y_n \uparrow \infty$ for $n \in \Z_+$ such that $2^{-n}G_1(y_n) \geq 
G_1(2 y_n)$ for $n \in \Z_+$. Take $x_n \uparrow \infty$ satisfying $x_0=0$ and \begin{eqnarray*}
 1 \leq G_1(y_n)(x_{n+1}-x_n) \leq 2 \quad\mbox{ for } n \in \Z_+.
\end{eqnarray*}
Furthermore, we define $g(x) \in {\cal G}_1 $ by $g(x) = y_n$ on $[x_n,x_{n+1})$. Then we obtain that
\begin{eqnarray*}
&& \int_0^{\infty}P(|X(1)-X(e^{-1})|  >g(x) )dx = \sum_{n=0}^{\infty}G_1(y_n)(x_{n+1}-x_n) = \infty 
\end{eqnarray*}
and
\begin{eqnarray*}
&& \int_0^{\infty}P(|X(1)-X(e^{-1})| > 2 g(x) )dx \le \sum_{n=0}^{\infty}2^{-n}G_1(y_n)(x_{n+1}-x_n) < \infty. \nonumber
\end{eqnarray*}
It follows from Proposition \ref{prop4.1} and Theorem \ref{th2.1}  that there is $C_0 \in [1,2(1-e^{-H})^{-1}]$ such that \eqref{eq:1.2} holds with $C=C_0$. Thus \eqref{eq:1.2} holds with $C=1$ by replacing $g(x)$ with $C_0g(x)$.   \qed
\medskip

\noindent
{\it Proof of Theorem \ref{th2.3}:}  By Theorem \ref{th3.2}, we have $F(r)\asymp G_1(r)$. As $G_1(r)\in OR$, the theorem holds from  Proposition \ref{prop4.1}.   \qed
\medskip

\noindent
{\it Proof of Corollary \ref{cor2.1}:}  The proof is obvious from Theorem \ref{th2.3}. \qed
\medskip

\noindent
{\it Proof of Corollary \ref{cor2.2}:} By using Proposition \ref{prop3.3}, we can obtain the corollary from Theorems \ref{th2.2} and  \ref{th2.3}. \qed

\noindent
{\it Proof of Example \ref{ex2.1}:} Use Corollary \ref{cor2.2} (ii) for
\begin{eqnarray*}
M(r) = \Gamma\left(\frac{m+1}{2}\right)\left(\sqrt{\pi}\Gamma\left(\frac{m}{2}\right)\right)^{-1}(1+r^2)^{-(m+1)/2}.
\end{eqnarray*}
Remaining assertion  is clear from the first assertion. \qed
\medskip

\noindent
{\it Proof of Theorem \ref{th2.4}:} By \eqref{eq:4.27}, we have
\begin{eqnarray}
\label{eq:4.29}
 &&\int_0^{\infty}G_1(\delta g(x))dx = \int_{\R^d}g^{-1}(\delta^{-1}|x|)\rho_1(dx)\quad\mbox{for $\delta>0$}. 
\end{eqnarray}
Suppose that $g^{-1}(x)+\log x \in OR$ and there is $\{X(t)\}$ such
that \eqref{eq:1.2} holds with $C=1$. Note that $\int_{|x|>1}\log |x|
\rho_1(dx)<\infty$. By the same way as in the proof of Theorem
\ref{th2.2}, we see from \eqref{eq:4.29}  that absurdity occurs. Thus if  $g^{-1}(x)+\log x \in OR$, then there is  no $\{X(t)\}$ such that \eqref{eq:1.2} holds with $C=1$.  

Conversely, we suppose that $g^{-1}(|x|)$ is submultiplicative on $\R^d$ and $g^{-1}(x) +\log x\notin OR.$ 
There is $x_n \uparrow \infty$ for $n \in \Z_+$ such that $x_0 =1$, $e^{2H} x_n <  x_{n+1}$, and $2^{-n}(g^{-1}(x_n)+\log x_n) \geq g^{-1}(e^{-H}x_n)+\log(e^{-H}x_n)$ for $n \in \Z_+$ and that, for $C_n : = \int_{x_n}^{e^H x_n}(g^{-1}(x)+\log x)x^{-1}dx$ with $n \in\Z_+,$
it holds that $\sum_{n=0}^{\infty}C_n^{-1} < \infty.$
Define $K(r)$ as
\begin{eqnarray*}
&& K(r) = 
\left\{
\begin{array}{ll}
\sum_{j=0}^{\infty} C_j^{-1}& \q  \mbox{ for } 0< r < e^Hx_0,\\ 
\sum_{j=n}^{\infty} C_j^{-1}& \q  \mbox{ for } e^Hx_{n-1} \le r < e^Hx_n \mbox{ with } n \in \N.\\ 
\end{array}
\right.
\end{eqnarray*}
Then we have
\begin{eqnarray*}
&& K(r)-K(e^Hr) = 
\left\{
\begin{array}{ll}
 C_n^{-1}& \q  \mbox{ for } x_n \le r < e^Hx_n \mbox{ with } n \in \Z_+,\\ 
0& \q \mbox{ for }  e^Hx_n \le r < x_{n+1} \mbox{ with } n \in \Z_+.\\ 
\end{array}
\right.
\end{eqnarray*}
Then we obtain that
\begin{eqnarray*}
&& \int_{1}^{\infty}(g^{-1}(r)+\log r)\frac{ K(r)-K(e^Hr)}{r}dr = \sum_{n=0}^{\infty}C_nC_n^{-1} = \infty
\end{eqnarray*}
and
\begin{eqnarray*}
& &\int_{1}^{\infty}(g^{-1}(e^{-2H}r)+\log(e^{-2H}r))\frac{ K(r)-K(e^Hr)}{r}dr \\
 &\le& \sum_{n=0}^{\infty}(g^{-1}(e^{-H}x_n)+\log(e^{-H}x_n))\int_{x_n}^{e^Hx_n}\frac{ K(r)-K(e^Hr)}{r}dr\nonumber  \\
 &\le& \sum_{n=0}^{\infty}2^{-n}\int_{x_n}^{e^Hx_n}(g^{-1}(r)+\log r)\frac{ K(r)-K(e^Hr)}{r}dr\nonumber  \\
 &= &\sum_{n=0}^{\infty}2^{-n}C_nC_n^{-1} < \infty.\nonumber 
\end{eqnarray*}
Notice that $\int_{1}^{\infty}(\log r)( K(r)-K(e^Hr))r^{-1} dr < \infty$ and $\int_{1}^{\infty}K(r)r^{-1} dr < \infty$. It follows from Theorem \ref{th2.1} (iii)  that \eqref{eq:1.2} holds with $C=C_a$ for some $C_a \in [1,e^{2H}]$. Hence \eqref{eq:1.2} holds with $C=1$ by replacing $\{X(t)\}$ with $\{C_a^{-1}X(t)\}$. \qed

\medskip

\noindent
{\it Proof of Theorem \ref{th2.5}:}  Note that $\int_{|x|>1}\log |x| \rho_1(dx)<\infty$. By using \eqref{eq:4.29}, the proof is clear from Proposition \ref{prop4.1}.   \qed

\medskip

\noindent
{\it Proof of Proposition \ref{prop2.1}:} Let $C\in [0,\infty).$ As mentioned in Remark \ref{rem3.1} (ii), $g^{-1}(|x|)$ is quasi-submultiplicative on $\R^d$. Thus we see from Proposition \ref{prop3.1} that \eqref{eq:2.23} holds if and only if \eqref{eq:2.13} holds. Here we assumed that $\int_{\R^d}g^{-1}(\delta^{-1}|x|)\mu(dx) < \infty$ for some $\delta \in (0,\infty)$ in the ``if'' part. Hence  assertion (i) holds  from Theorem \ref{th2.1}. Now we see from Lemma \ref{lem3.2} that \eqref{eq:2.15} is equivalent to \eqref{eq:2.23}. Furthermore, we see that $\int_{1}^{\infty} K(r)g^{-1}(\delta^{-1}r)r^{-1} dr < \infty$ if and only if $\int_{\R^d}g^{-1}(\delta^{-1}|x|)\mu(dx) < \infty$. Thus assertion (ii) holds from (i). \qed
\medskip

\noindent
{\it Proof of Proposition \ref{prop2.2}:} Define $\zeta_1$ on $\R_+$ by $\zeta_1(x >r) := \rho_1(|x|>r)$.
Suppose that $\rho_1 \in {\mathcal O}{\mathcal S}$ on $\R^d$, that is, $\zeta_1 \in {\mathcal O}{\mathcal S}$ on $\R_+$. Now we have, for $g \in {\cal G}_1 $  and $l\in \N$, 
\begin{eqnarray*}
 \int_{\R^d}g^{-1}(|x|)\rho_l(dx)& =& \int_{(\R^d)^l}g^{-1}\left(\left|\sum_{k=1}^le^{(-k+1)H}x_k\right|\right)\prod_{k=1}^l\rho_1(dx_k). \\
& \le &\int_{(\R^d)^l}g^{-1}(\sum_{k=1}^l|x_k|)\prod_{k=1}^l\rho_1(dx_k)\nonumber\\
& = &\int_{\R_+}g^{-1}(r)\zeta_1^{l*}(dr)\nonumber\\
& = &\int_0^{\infty}\zeta_1^{l*}(x >g(r))dr.\nonumber
\end{eqnarray*}
By Lemma \ref{lem5.5} (iii) in Sect.5, $\zeta_1^{l*}(x >r) \asymp
\zeta_1(x >r)$. Thus, if $\int_{\R^d}g^{-1}(|x|)\rho_1(dx) < \infty$,
then $\int_{\R^d}g^{-1}(|x|)\rho_l(dx)< \infty$ for all $l\in
\N$. Therefore we see that \eqref{eq:2.13} holds if and only if
\eqref{eq:2.25} holds. By virtue of Theorem \ref{th3.1} (iv),
assertion (ii) holds from (i). The proof of (iii) is clear from
Proposition \ref{prop2.1} (i) and Theorem \ref{th3.1} (iv). \qed

\medskip

\noindent
{\it Proof of Example \ref{ex2.2}:}  The example of Sect.6 in \cite{egv79} shows that $\mu \in {\mathcal S}$ on $\R_+$ and 
\begin{eqnarray}
\label{eq:4.30}
\mu(x >r) \sim (2\pi)^{-1/2}(\log r)^{-1}\exp(-(\log r)^2/2).
\end{eqnarray}
Hence we see from Theorem 1 of \cite{egv79} that
\begin{eqnarray*}
\int_r^{\infty}\frac{K(x)}{x}dx \sim \mu(x >r).
\end{eqnarray*}
By \eqref{eq:2.11} and \eqref{eq:4.30}, we have 
\begin{eqnarray}
\label{eq:4.31}
\eta_1(x >r) = \int_r^{\infty}\frac{K(x)}{x}dx - \int_{e^Hr}^{\infty}\frac{K(x)}{x}dx \sim  \mu(x >r).
\end{eqnarray}
It follows from Lemma A3.15 of \cite{ekm97} that $\widetilde\eta_1 \in {\mathcal S}$, and thereby  $\widetilde\eta_1 \in {\mathcal O}{\mathcal S}$. 
For $\delta >0$ and $g \in{\mathcal G}_1$, we find from \eqref{eq:4.31} that 
$$\int_{\R_+}\eta_1(x>\delta g(r))dr=\int_{\R_+}g^{-1}(\delta^{-1} x)\eta_1(dx) < \infty$$
 if and only if
\begin{eqnarray*}
\int_{\R_+}\mu(x>\delta g(r))dr=\int_{\R_+}g^{-1}(\delta^{-1} x)\mu(dx) < \infty.
\end{eqnarray*}
Thus we obtain the first assertion from Proposition \ref{prop2.2} (ii). Setting $g^{-1}(x) = \exp((\log x)^2/2)$ for $x>1$, we obtain the second assertion from the first one. The second assertion can be proved also by employing Proposition \ref{prop2.1} (i).\qed

\medskip

\begin{lem} 
\label{lem4.8}
Let $C \in [0,\infty]$. 
 Let $\phi(r)$ be a positive, increasing, and regularly  varying function with positive index. Suppose that $g \in {\cal G}_1$ satisfies that
\begin{eqnarray}
&& g^{-1}(r) = \exp(\phi(r)).
\end{eqnarray}
Then \eqref{eq:2.13} holds if and only if
\begin{eqnarray}
\label{eq:4.33}
&& \int_1^{\infty}\rho_l(|x| >r) g^{-1}(\delta^{-1}r)dr \nonumber\\
&&\quad \left\{\begin{array}{l}
<\infty \quad  \mbox{ for $\delta>C$ and all $l\in {\N}$},\\
=\infty \quad  \mbox{ for $\delta<C$ and  some $l= l(\delta)\in {\N}$}.
\end{array}
\right.
\end{eqnarray}
Moreover, if $g^{-1}(r)$ is submultiplicative on $\R_+$, then \eqref{eq:2.13} is also equivalent to \eqref{eq:2.23}, and to 
\begin{eqnarray}
\label{eq:4.34}
&& \int_1^{\infty}K(r) g^{-1}(\delta^{-1}r)dr 
\left\{\begin{array}{l}
<\infty \quad  \mbox{ for }    \delta>C,\\
=\infty \quad   \mbox{ for }     0< \delta<C.
\end{array}
\right.
\end{eqnarray}
\end{lem}

\prf By Theorem 4.12.10 (ii) of \cite{bgt87}, \eqref{eq:2.13} holds if and only if  
\begin{eqnarray}
\label{eq:4.35}
&&\int_{\R^d}\rho_l(dx) \int_0^{|x|}g^{-1}(\delta^{-1}r)dr \nonumber\\
&&\quad \left\{\begin{array}{l}
<\infty \quad   \mbox{ for $\delta>C$ and all $l\in {\N}$},\\
=\infty \quad   \mbox{ for $0<\delta<C$ and some $l=l(\delta)\in {\N}$}.\qquad
\end{array}
\right.
\end{eqnarray}
Here we used the assumption that $\phi(r)$ is regularly varying. The integral of \eqref{eq:4.35} is equal to that of \eqref{eq:4.33}. \par
Suppose that $g^{-1}(r)$ is submultiplicative on $\R_+$. Then $\int_0^{|x|\vee 1}g^{-1}(r)dr$ is submultiplicative on $\R^d$ and thus, by Lemma \ref{lem3.2}, \eqref{eq:4.35} is equivalent to
\begin{eqnarray*}
&& \int_{\R^d}\eta_l(dx) \int_0^{|x|}g^{-1}(\delta^{-1}r)dr \\
&&\quad \left\{\begin{array}{l}
<\infty \quad  \mbox{ for $\delta>C$ and all $l\in {\N}$},\\
=\infty \quad  \mbox{ for $0<\delta<C$ and some $l=l(\delta)\in {\N}$},
\end{array}
\right.
\end{eqnarray*}
equivalently,
 \begin{eqnarray}
\label{eq:4.36}
&& \int_1^{\infty}\eta_l(|x| >r) g^{-1}(\delta^{-1}r)dr \nonumber\\
&&\quad \left\{\begin{array}{l}
<\infty \quad \mbox{ for $\delta>C$ and all $l\in {\N}$},\\
=\infty \quad \mbox{ for $0<\delta<C$ and some $l=l(\delta)\in {\N}$}.
\end{array}
\right.\qquad
\end{eqnarray}
By \eqref{eq:2.11}, we have
 \begin{eqnarray}
\label{eq:4.37}
&& lH K(r) \ge \eta_l(|x| >r) \ge (\log(1+\ep)) K((1+\ep)r)
\end{eqnarray}
for sufficiently small  $\ep >0$. Thus \eqref{eq:4.36} is equivalent to \eqref{eq:4.34}. Moreover, by Lemma \ref{lem3.2}, \eqref{eq:2.23} is equivalent to \eqref{eq:2.15}. As $\phi(x)$ is regularly varying  with positive index, \eqref{eq:2.15} is equivalent to \eqref{eq:4.34}. Thus the second assertion is true. \qed

\medskip

\noindent
{\it Proof of Proposition \ref{prop2.3}:} Let 
\begin{eqnarray*}
g^{-1}(r) = \exp(r^{\alpha}f(\log(r\vee 1))).
\end{eqnarray*}
Since $f(r)$ is regularly varying, $f(\log(r\vee 1))$ is slowly varying. Thus we have
\begin{eqnarray*}
g(r) \sim \frac{ (\log r)^{1/\alpha}}{(f(\alpha^{-1}\log_{(2)} r))^{1/\alpha}}.
\end{eqnarray*}
Notice that \eqref{eq:2.27} and \eqref{eq:2.28} are equivalent  to \eqref{eq:4.34} and \eqref{eq:2.23}, respectively.
Since $g^{-1}(r)$ is submultiplicative on $\R_+$, we see from Lemma \ref{lem4.8} that \eqref{eq:2.26} is equivalent to \eqref{eq:2.27} and to \eqref{eq:2.28}. \qed

\medskip

\noindent
{\it Proof of Proposition \ref{prop2.4}:}  Let 
\begin{eqnarray*}
g^{-1}(r) = \exp(rf^{-1}(\log (r\vee 1))).
\end{eqnarray*}
We see  from Definition \ref{def2.4} that, for all sufficiently small $\ep >0$,
\begin{eqnarray*}
\limsup_{x \to \infty}\frac{f^{-1}((1+\ep)x)}{ f^{-1}(x)} \le 1+\ep.
\end{eqnarray*}
Thus  $f^{-1}(\phi(t)) \sim f^{-1}(t)$ provided that $\phi(t) \sim t$ as $t \to \infty.$ This implies that $rf^{-1}(\log (r+1))$ is regularly varying with index 1, and we have
\begin{eqnarray*}
g(r) \sim   \frac{ \log r}{f^{-1}(\log_{(2)} r)}.
\end{eqnarray*}
By Theorem \ref{th2.1} and Lemma \ref{lem4.8}, \eqref{eq:2.29} is equivalent to \eqref{eq:4.33}, that is,
\begin{eqnarray*}
& & \int_1^\infty\exp(\delta^{-1} rf^{-1}(\log (r+1)))\rho_l(|x| >r)dr\\
 &   &\left\{
\begin{array}{ll}
<\infty \quad   \forall l\in {\N},  \mbox{ for } \delta>C,\nonumber\\
=\infty \quad   \exists l=l(\delta)\in {\N},   \mbox{ for }0<\delta<C.\nonumber
\end{array}
\right.
\end{eqnarray*}
By Lemma \ref{lem3.4}, this is equivalent to
\begin{eqnarray}
\label{eq:4.38}
 &  &\int_1^\infty\exp(h(\delta^{-1} r))\eta_l(|x| >r) dr\nonumber\\
 &  & \left\{
\begin{array}{ll}
<\infty \quad   \forall l\in {\N},  \mbox{ for } \delta>C,\\
=\infty \quad   \exists l=l(\delta)\in {\N},   \mbox{ for }0<\delta<C.
\end{array}
\right.
\end{eqnarray}
Using \eqref{eq:4.37} for sufficiently small $\ep>0$, we see that
\eqref{eq:4.38} is equivalent to
\begin{eqnarray}
\label{eq:4.39}
\int_1^\infty\exp(h(\delta^{-1} r))K(r)dr
 \left\{
\begin{array}{ll}
<\infty \quad   \mbox{ for } \delta>C,\\
=\infty \quad    \mbox{ for }0<\delta<C.
\end{array}
\right.
\end{eqnarray}
Using \eqref{eq:4.37} again, we see that this is equivalent to
\begin{eqnarray}
\label{eq:4.40}
   \int_1^\infty\exp(h(\delta^{-1} r))\eta_1(|x| >r) dr
    \left\{
\begin{array}{ll}
<\infty \quad     \mbox{ for } \delta>C,\\
=\infty \quad      \mbox{ for }0<\delta<C.
\end{array}
\right.
\end{eqnarray}
Further, \eqref{eq:4.39} is equivalent to \eqref{eq:2.30}. We see from Lemma \ref{lem3.4} that \eqref{eq:4.40} is equivalent to \eqref{eq:2.31}.

Lastly, we prove the equivalence of \eqref{eq:2.30} and \eqref{eq:2.32} provided that $h(r)$ is regularly varying with positive index. Now  \eqref{eq:2.30} is equivalent to \eqref{eq:4.39}, equivalently,
\begin{eqnarray*}
\int_1^\infty\exp(h(\delta^{-1} r))K(r)r^{-1}dr
 \left\{
\begin{array}{ll}
<\infty \quad    \mbox{ for } \delta>C,\\
=\infty \quad   \mbox{ for }0<\delta<C.
\end{array}
\right.
\end{eqnarray*}
Here we used the regularly variation. By virtue of Theorem 4.12.10 (ii) of \cite{bgt87}, this is equivalent to
\begin{eqnarray*}
   \int_1^\infty\left(\int_0^r\exp(h(\delta^{-1} s))ds\right)K(r)r^{-1}dr
   \left\{
\begin{array}{ll}
<\infty \quad   \mbox{ for } \delta>C,\\
=\infty \quad    \mbox{ for }0<\delta<C,
\end{array}
\right.
\end{eqnarray*}
 equivalently,
\begin{eqnarray}
\label{eq:4.41}
   \int_1^\infty \left( \int_r^{\infty}K(s)s^{-1}ds\right) \exp(h(\delta^{-1} r))dr
    \left\{
\begin{array}{ll}
<\infty \quad   \mbox{ for } \delta>C,\\
=\infty \quad    \mbox{ for }0<\delta<C.
\end{array}
\right.
\end{eqnarray}
Notice that $\nu(|x|>r)=\int_r^\infty K(s)s^{-1}ds$. It follows that \eqref{eq:4.41} is equivalent to
$$\liminf_{r\to\infty}\frac{h^{-1}(-\log \nu(|x|>r))}{r}=C.$$
By Lemma \ref{lem3.4}, this is equivalent to \eqref{eq:2.32}.  \qed

\medskip

\noindent
{\it Proof of Proposition \ref{prop2.5}:} Let $g^{-1}(r) = \exp(r\log (r\vee 1))$. The proof of  Proposition \ref{prop2.5} is similar to that of  Proposition \ref{prop2.4} by using Lemma \ref{lem3.3} in place of Lemma \ref{lem3.4}. It is omitted.\qed
\medskip

\noindent
{\it Proof of Theorem \ref{th2.6}:}
Note from Theorem \ref{th2.1} (ii) that if $g_1^{-1}(r) \asymp  g_2^{-1}(r)$ for $g_1, g_2 \in {\cal G}_1,$ then \eqref{eq:1.2} holds for $g=g_1$ and for $g=g_2$ with the same constant $C$.  Thus in assertion (i), (1), (2), and (3) are from Remark \ref{rem2.5}, Proposition \ref{prop2.5} and Proposition \ref{prop2.4}, respectively.\par

 Next we prove assertion (ii). Suppose that \eqref{Th2.6.1} holds. Since $\rho_1$ is infinitely divisible, we see from Lemma \ref{lem3.3} that 
\begin{eqnarray*}
\int_{\R^d} e^{c|x|^2}\rho_1(dx) = \infty \q \mbox{ for  some }c >0.
\end{eqnarray*}
Hence, by \eqref{Th2.6.1}, we have 
\begin{eqnarray}
\int_{\R^d} g^{-1}(\delta^{-1} |x|)\rho_1(dx)  = \infty \q \mbox{ for  } 0 <\delta < \infty. \label{PTh2.6.1}
\end{eqnarray}
By Theorem \ref{th2.1} (ii),  \eqref{eq:1.2} holds with $C=\infty$. Hence assertion (1) is true. Suppose that \eqref{Th2.6.2} and \eqref{Th2.6.3} hold. In the case where $\nu \ne 0$, we find from Lemma \ref{lem3.3} that
\begin{eqnarray*}
\int_{|x| >1} e^{c|x|\log |x|}\rho_1(dx) = \infty \q \mbox{ for  some }c >0.
\end{eqnarray*}
Thus we have \eqref{PTh2.6.1} by \eqref{Th2.6.2}, and thereby \eqref{eq:1.2} holds with $C=\infty$. In the case where $\nu =0$, we see that
\begin{eqnarray*}
\int_{\R^d} e^{c|x|^2}\mu(dx) < \infty \mbox{ for  some } c >0.
\end{eqnarray*}
Hence, by \eqref{Th2.6.3}, we have
\begin{eqnarray*}
\int_{\R^d} g^{-1}(\delta^{-1} |x|)\mu(dx)  < \infty \q \mbox{ for  } 0 <\delta < \infty.
\end{eqnarray*}
Hence we obtain from Proposition \ref{prop2.1} that \eqref{eq:1.2} holds with $C=0$. Hence assertion (2) is true.\qed
\medskip

\noindent
{\it Proof of Example \ref{ex2.3}:} (i) Use Proposition \ref{prop2.3} with $f(x) = x^{\beta}$.\\
(ii) Use Proposition \ref{prop2.4} with $C= c^{-1/\alpha}$ and $h(x) = x^{\alpha}(\log x)^{\beta}$ for sufficiently large $x$. Then it suffices to prove that
\begin{eqnarray*}
&& f^{-1} (y) \sim  
 \left\{
\begin{array}{ll}
(\log y)^{\beta}\q  \mbox{ for }\alpha = 1,\\
 \alpha^{-(\beta-\alpha)/\alpha}(\alpha-1)^{-(\alpha-1)/\alpha} y^{(\alpha-1)/\alpha}(\log y)^{\beta/\alpha}\q  \mbox{ for }\alpha > 1. 
\end{array}
\right.
\end{eqnarray*}
Suppose that $\alpha>1$. As $h(x)=\int_0^x q(t)dt$, it follows that $q(t)\sim \alpha t^{\alpha-1}(\log t)^\beta$ as $t\to\infty$. Put $y=q(t)$. Then $y\sim \alpha t^{\alpha-1}(\log t)^\beta$ and $\log y\!\sim\! (\alpha-1)\log t$ as $t\to\infty$. Hence
\begin{eqnarray*}
t\sim\left(\frac{y}{\alpha (\log t)^\beta}\right)^{1/(\alpha-1)}\sim\left(\frac{y}{\alpha ((\alpha-1)^{-1}\log y)^\beta}\right)^{1/(\alpha-1)}
\end{eqnarray*}
as $t\to\infty$. This implies that 
\[q^{-1}(y)\!\sim\!
\alpha^{-1/(\alpha-1)}(\alpha-1)^{\beta/(\alpha-1)}y^{1/(\alpha-1)}(\log
y)^{-\beta/(\alpha-1)}.
\] 
Hence,
\begin{eqnarray*}
f(x)=\int_0^x q^{-1}(y)dy &\sim& \alpha^{-1/(\alpha-1)}(\alpha-1)^{\beta/(\alpha-1)}\frac{x^{(\alpha-1)^{-1}+1}}{(\alpha-1)^{-1}+1}(\log x)^{-\beta/(\alpha-1)}\\
&=& \alpha^{-\alpha/(\alpha-1)}(\alpha-1)^{(\alpha+\beta-1)/(\alpha-1)}x^{\alpha/(\alpha-1)}(\log x)^{-\beta/(\alpha-1)}.
\end{eqnarray*}
Moreover, put $y=f(x)$. Then we have $\log y\sim \alpha(\alpha-1)^{-1}\log x$ as $x\to\infty$. Hence,
\begin{eqnarray*}
x &\sim& \left(\frac{\alpha^{\alpha/(\alpha-1)}(\alpha-1)^{-(\alpha+\beta-1)/(\alpha-1)}y}{(\log x)^{-\beta/(\alpha-1)}}\right)^{(\alpha-1)/\alpha}\\
&\sim& \left(\frac{\alpha^{\alpha/(\alpha-1)}(\alpha-1)^{-(\alpha+\beta-1)/(\alpha-1)}y}{((\alpha-1)\alpha^{-1}\log y)^{-\beta/(\alpha-1)}}\right)^{(\alpha-1)/\alpha}\\
&=& \alpha^{-(\beta-\alpha)/\alpha}(\alpha-1)^{-(\alpha-1)/\alpha}y^{(\alpha-1)/\alpha}(\log y)^{\beta/\alpha}.
\end{eqnarray*}
This implies that $f^{-1}(y)\sim\alpha^{-(\beta-\alpha)/\alpha}(\alpha-1)^{-(\alpha-1)/\alpha}y^{(\alpha-1)/\alpha}(\log y)^{\beta/\alpha}$. Next suppose that $\alpha=1$. Then $y=q(t)=(\log t)^\beta+\beta(\log t)^{\beta-1}$ for all sufficiently large $t$. Hence,
\begin{eqnarray*}
\exp(y^{\beta^{-1}})=\exp((\log t)(1+\beta(\log t)^{-1})^{\beta^{-1}})\sim e\cdot t\quad\mbox{as $t\to\infty$.}
\end{eqnarray*}
This yields that $q^{-1}(y)\sim e^{-1}\exp(y^{\beta^{-1}})$. Hence we obtain that
\begin{eqnarray*}
f(x)=\int_0^x q^{-1}(y)dy \sim e^{-1}\exp(x^{\beta^{-1}})\cdot \beta x^{1-\beta^{-1}}.
\end{eqnarray*}
Put $y=f(x)$. Then $\log y\sim x^{\beta^{-1}}$. Hence we obtain that $f^{-1}(y)\sim(\log y)^\beta$. Thus all assertions are true. \qed


\section{Proof of Theorem \ref{th3.1}.}

In this section, we prove Theorem \ref{th3.1} mentioned in Sect.3. The proof of Theorem \ref{th3.1} is similar to that of Theorem 1.1 of \cite{sw05}, but there is a difficulty peculiar to the multi-dimensional case.

\begin{lem}
\label{lem5.1}
 Let $\rho_j$ for $j=1,2$ be distributions on $\R^d$. For $j=1,2$, we define a distribution $\zeta_j$ on $\R_+$ by $\zeta_j(x>r):=\rho_j(|x|>r)$ for $r \ge 0$.

{\rm (i)} We have 
\begin{equation}
\rho_1*\rho_2(|x| >r) \le \zeta_1*\zeta_2(x>r)\q \mbox{ for }  r \ge 0.
\end{equation}
In particular, suppose $n \in \N $, then we have
\begin{equation}
\rho_1^{n*}(|x| >r) \le \zeta_1^{n*}(x>r)\q\mbox{ for } r \ge 0.
\end{equation}

{\rm (ii)}  There are  $ s >0$ and $c_1 >1$ both independent of  $\rho_1$  such that      
\begin{equation}
\rho_1(|x| >r) \le c_1 \rho_1*\rho_2(|x| > r-s)\q \mbox{ for }  r \ge 0.
\end{equation}
\end{lem}

\prf Let $\{X_j\}$ be independent $\R^d$-valued random variables such that  the  distribution of $X_j$ is $\rho_j$ for $j = 1,2.$ Then $\zeta_j$ is the  distribution of $|X_j|$ for $j = 1,2.$ Thus we have, for  $r \ge 0$,
\begin{eqnarray*}
\rho_1*\rho_2(|x| >r)& =& P(|X_1+X_2| >r)\\
  & \le&  P(|X_1|+|X_2| >r)=\zeta_1*\zeta_2(x>r).\nonumber
\end{eqnarray*}
The second assertion of (i) is trivial.

Choose $s >0$ such that $c_1^{-1}:= P(|X_2| \le s) >0.$
Then we see that for $r >s$,
\begin{eqnarray*}
\rho_1(|x| >r) & =& c_1P(|X_2| \le s)P(|X_1| >r)\\
& =& c_1P(|X_1| >r,|X_2| \le s)\nonumber\\
& \le& c_1  P(|X_1+X_2| >r-s)=c_1 \rho_1*\rho_2(|x| >r-s).\nonumber
\end{eqnarray*}
We have proved (ii).\qed

\begin{lem}
\label{lem5.2}
 Let $\zeta_j$ for $j=1,2,3$ be distributions on $\R_+$.
 If $\zeta_1(x >r ) \le c_1 \zeta_2(x >r ) $ for some $c_1 >1,$
then 
\begin{equation}
\label{eq:5.4}
\zeta_1*\zeta_3(x >r ) \le c_1  \zeta_2*\zeta_3(x >r)\q \mbox{ for }  r \ge 0,
\end{equation}
and, for any $n\in \N$, 
\begin{equation}
\label{eq:5.5}
\zeta_1^{n*}(x >r) \le c_1^n \zeta_2^{n*}(x >r)\q \mbox{ for }  r \ge 0.
\end{equation}
\end{lem}

\prf Suppose that $\zeta_1(x >r ) \le c_1 \zeta_2(x >r ) $ for some $c_1 >1.$ Then we see
\begin{eqnarray*}
\zeta_1*\zeta_3(x >r )  & =& \int_{0-}^{r+}\zeta_1(x >r-u )\zeta_3(du)+\zeta_3(x >r )  \\
& \le&  c_1\int_{0-}^{r+}\zeta_2(x >r-u )\zeta_3(du) + c_1\zeta_3(x >r ) \nonumber\\
& =& c_1 \zeta_2*\zeta_3(x >r).\nonumber
\end{eqnarray*}
The inequality \eqref{eq:5.5} is trivial from \eqref{eq:5.4}. \qed

\begin{lem}
\label{lem5.3}
 Let $\zeta_j$ for $j=1,2$ be distributions on $\R_+$. 

{\rm (i)} If $\zeta_1 \in {\mcal O}{\mcal S}$  and $\zeta_1(x >r) \asymp \zeta_2(x >r)$, then  $\zeta_2\in {\mcal O}{\mcal S}$.

{\rm (ii)} If $\zeta_1 \in {\mcal D}$  and $\zeta_1(x >r) \asymp \zeta_2(x >r)$, then  $\zeta_2\in {\mcal D}$.
\end{lem}

\prf Assertion (i) is from Theorem 2.3 of \cite{kl88}. Assertion (ii) is clear from the definition. \qed

\begin{lem}
\label{lem5.4}
 Let $\zeta$  be distribution on $\R_+$. Then $\zeta^{n*} \in {\mcal D}$ for some $n \in \N$ if and only if  $\zeta \in {\mcal D}$.
\end{lem}

\prf This is from Proposition 1.1 (iii) and Proposition 2.5 (iii) of \cite{sw05}. \qed
\medskip

For  $\zeta \in {\mathcal O}{\mathcal S}$ on $\R_+$, we define 
\begin{equation}
\ell^*(\zeta):=\limsup_{r \to \infty}\frac{\zeta^{2*}(x>r)}{\zeta(x>r)} < \infty.
\end{equation}

\medskip

\begin{lem}
\label{lem5.5}
 Suppose that $\zeta \in {\mcal O}{\mcal S}$ on $\R_+$.

{\rm (i)} For any $s >0$, there is $c_1 >0$ such that $\zeta(x >r) \le c_1\zeta(x >r+s)$. Thus $\zeta(x >\log(r+1)) \in D$ and there are $c_2,c_3 >0$ such that  
\begin{equation}
\zeta(x >r)\ge c_2\exp(-c_3 r)\q \mbox{ for } r \ge 0.
\end{equation}

{\rm (ii)} For every $\ep >0$ there is $c_4 >0$ such that
\begin{equation}
\label{eq:5.8}
\zeta^{k*}(x>r)  \le c_4(\ell^*(\zeta)-1 +\ep)^k\zeta(x>r)\mbox{ for every  } r \ge 0  \mbox{ and } k\ge1.
\end{equation}

{\rm (iii)} For every $k \in \N$, we have
\begin{equation}
\zeta^{k*}(x>r) \asymp \zeta(x>r).
\end{equation}
\end{lem}

\prf Assertion (i) is from Proposition 2.2 (ii) of \cite{sw05}. Assertion (ii) is from  Proposition 2.4  of \cite{sw05}. As $\zeta^{k*}(x>r) \ge \zeta(x>r)$ for every $k \in \N$, assertion (iii) is clear from (ii). \qed

\medskip

\begin{lem}
\label{lem5.6}
 Let $\rho$ be an infinitely divisible distribution on  ${\R}^d$ with L\'evy measure $\eta$. Take $c>0$ such that $\eta(|x|>c)>0$. Let  $\rho_j$ for $j=1,2$ be the  infinitely divisible distributions on  ${\R}^d$ such that $\rho_2$ is a compound Poisson distribution on ${\R}^d$ with  L\'evy measure $1_{\{|x| >c\}}\eta(dx)$ and $\rho= \rho_1*\rho_2$. Then the following hold:

{\rm (i)}  $\rho\in {\mathcal O}{\mathcal S}$ if and only if $\rho_2\in {\mathcal O}{\mathcal S}$.

{\rm (ii)} If  $\rho\in {\mathcal O}{\mathcal S}$, then $\rho(|x| > r) \asymp \rho_2(|x| > r)$.
\end{lem}

\prf Define  the distributions $\zeta_j$ for $j = 1,2$ and
$\zeta_{\rho}$ by $\zeta_j(x >r )= \rho_j(|x| >r ) $ and
$\zeta_{\rho}(x >r) = \rho(|x| >r)$ for $r \ge 0$. Note from Lemma
\ref{lem3.3} that $\zeta_1(x >r ) \le c_1 \zeta_2(x >r )$ for some
$c_1 >1$. On the one hand, we see from Lemmas  \ref{lem5.1} (i) and
\ref{lem5.2} that 
\begin{equation}
\rho(|x| >r) = \rho_1*\rho_2(|x| >r) \le \zeta_1*\zeta_2(x>r)\le c_1 \zeta_2^{2*}(x>r).\label{lem5.6.1}
\end{equation}
On the other hand, we see from Lemma \ref{lem5.1} (ii)  that, for some $c_2 \in (0,1)$ and $s>0,$  
\begin{equation}
\rho(|x| >r) = \rho_1*\rho_2(|x| >r) \ge c_2 \zeta_2(x>r+s).\label{lem5.6.2}
\end{equation}
Suppose that $\rho_2\in {\mathcal O}{\mathcal S}$, that is, $\zeta_2\in {\mathcal O}{\mathcal S}$. Thus 
$\zeta_2^{2*}(x>r)\asymp \zeta_2(x >r ).$ 
Then, by using Lemma \ref{lem5.5} (i)  for $\zeta_2$,   we obtain from \eqref{lem5.6.1} and \eqref{lem5.6.2} that $\rho(|x| > r) \asymp \zeta_2(x > r),$ and, by Lemma \ref{lem5.3} (i),  $\rho\in {\mathcal O}{\mathcal S}$. Conversely, suppose that  $\rho\in {\mathcal O}{\mathcal S}$. Then $\zeta_{\rho}^{2*}(x>r)\le  c_3\zeta_{\rho}(x >r )$ for some $ c_3 >0$. By using Lemma \ref{lem5.5} (i) for $\zeta_{\rho}$, we have, for some $ c_4 >0$,
\begin{equation}
\rho(|x| >r) \ge c_4\rho(|x| >r-s) \ge c_2c_4 \zeta_2(x>r).\label{lem5.6.3}
\end{equation}
This implies that there is $ c_5 >1$ such that $\zeta_2(x>r) \le c_5 \zeta_{\rho}(x >r)$ for $r \ge0$. From Lemma \ref{lem3.3}, we can take $\ep >0$ and $A>0$ such that $\ep c_3 c_5  <1$ and $\zeta_1(x >r ) < \ep \zeta_{\rho}(x >r)$ for $r \ge A$. Hence we see that {\allowdisplaybreaks
\begin{eqnarray*}
\rho(|x| >r+A)& \le&\zeta_1*\zeta_2(x>r+A)\\
& \le&\int_{0-}^{r+}\zeta_1(x > r+A-u)\zeta_2(du) +\zeta_2(x>r)\nonumber\\
& \le&\ep \int_{0-}^{r+}\zeta_{\rho}(x > r+A-u)\zeta_2(du) +\zeta_2(x>r)\nonumber\\
&\le&\ep \zeta_{\rho}*\zeta_2(x>r+A) +\zeta_2(x>r)\nonumber\\
&\le&\ep c_5 \zeta_{\rho}*\zeta_{\rho}(x>r+A) +\zeta_2(x>r)\nonumber\\
&\le&\ep c_3 c_5 \zeta_{\rho}(x>r+A) +\zeta_2(x>r).
\end{eqnarray*}
This yields }
\begin{eqnarray*}
(1- \ep c_3 c_5 )\rho(|x| >r+A) \le \zeta_2(x>r).
 \end{eqnarray*}
We obtain from Lemma \ref{lem5.5} (i) and \eqref{lem5.6.3}  that $\rho(|x| > r) \asymp \rho_2(|x| > r)$, and consequently $\rho_2\in {\mathcal O}{\mathcal S}$. Here also assertion (ii) has been proved. \qed

\medskip

\begin{prop}
\label{prop5.1}
 Let $\rho$ be a compound Poisson distribution on ${\R}^d$ with L\'evy measure $\eta$.
Define a distribution $\zeta$ on $\R_+$ satisfying  $\zeta(x >r)= \widetilde\eta(|x| >r)$ for $r \ge 0.$

{\rm (i)} $\rho\in {\mathcal O}{\mathcal S}$  on ${\R}^d$ if and only if there is a positive integer $n$ such that $\widetilde\eta^{n*}\in {\mathcal O}{\mathcal S}$ on $\R^d$. Moreover, if $\rho\in {\mathcal O}{\mathcal S}$, then \eqref{eq:3.12} holds.

{\rm (ii)} If $\widetilde\eta  \in {\mathcal O}{\mathcal S}$, then \eqref{eq:3.13} holds.

{\rm (iii)}  $\rho\in {\mathcal O}{\mathcal S}$ if and only if there is a positive integer $m$ such that $\zeta^{m*}\in {\mathcal O}{\mathcal S}$ on $\R_+$ and \eqref{eq:3.14} holds.
\end{prop}

\prf We define   distributions  $\zeta_{\rho}$ and  $\zeta_n$ on $\R_+$ as follows:   For $r \ge 0$,
 \begin{equation*}
\zeta_{\rho}(x >r): = \rho(|x| >r)\quad\mbox{and}\quad \zeta_n(x >r):= \widetilde\eta^{n*}(|x| >r).
\end{equation*}

 Let $\delta:= \eta(\R^d) < \infty$. Suppose that there is  $n\in\N$ such that $\widetilde\eta^{n*}\in {\mathcal O}{\mathcal S}$ on $\R^d$, namely, $\zeta_n\in {\mathcal O}{\mathcal S}$ on $\R_+$. 
Then we see from Lemma \ref{lem5.5} (ii) that there is $c_1 >0$ such that
\begin{equation*}
\zeta_n^{k*}(x>r)  \le c_1\ell^*(\zeta_n)^{k}\zeta_n(x>r)\mbox{ for every  } r \ge 0  \mbox{ and } k\ge1.
\end{equation*}
By using Lemma \ref{lem5.1} (ii)  for $\widetilde\eta^{(nk+j)*}$ with $0 \le j \le n-1$ and $s >0$, we obtain from Lemmas  \ref{lem5.1} (i) and \ref{lem5.5} (i)  that, for some $c_2, c_3 >0$, {\allowdisplaybreaks
\begin{eqnarray*}
\rho(|x| >r)& =& e^{ -\delta}\sum_{k=0}^{\infty}\frac{\delta^k}{k!}\widetilde\eta^{k*}(|x|>r)\\
            & =& e^{ -\delta}\sum_{j=0}^{n-1}\sum_{k=0}^{\infty}\frac{\delta^{nk+j}}{(nk+j)!}\widetilde\eta^{(nk+j)*}(|x|>r)\nonumber\\
            &\le& c_2\sum_{k=0}^{\infty}\frac{\delta^{nk}}{(nk)!}\widetilde\eta^{n(k+1)*}(|x|>r-s)\nonumber\\
          &\le & c_2\sum_{k=0}^{\infty}\frac{\delta^{nk}}{(nk)!}\zeta_n^{(k+1)*}(x>r-s)\nonumber\\
            &\le & c_2c_1\sum_{k=0}^{\infty}\frac{\ell^*(\zeta_n)^{k+1}\delta^{nk}}{(nk)!}\zeta_n(x>r-s)\\
            &\le& c_3\zeta_n(x>r).\nonumber
 \end{eqnarray*}
Obviously, it follows that}
\begin{equation}
\rho(|x| >r)\ge  e^{ -\delta}\frac{\delta^n}{n!}\widetilde\eta^{n*}(|x|>r)=  e^{ -\delta}\frac{\delta^n}{n!}\zeta_n(x>r).\label{Prop5.1.1}
\end{equation}
Thus we have
\begin{equation}
\rho(|x| >r)\asymp  \zeta_n(x>r)=\widetilde\eta^{n*}(|x|>r).\label{Prop5.1.2}
\end{equation}
Hence $\rho \in {\mathcal O}{\mathcal S}$ by Lemma \ref{lem5.3} (i). 

Conversely, we suppose that $\rho \in {\mathcal O}{\mathcal S}$ on $\R^d$. By using the method in the proof of Theorem 1.5 of \cite{wa08}, we shall prove that there is  $n\in\N$ such that
\begin{equation}
\liminf_{r \to \infty}\frac{\widetilde\eta^{n*}(|x|>r)}{\rho(|x| >r)} >0.\label{Prop5.1.3}
\end{equation}
Suppose that, for all $n\in\N$,  
\begin{equation*}
\liminf_{r \to \infty}\frac{\widetilde\eta^{n*}(|x|>r)}{\rho(|x| >r)} =0.
\end{equation*}
Then we can choose $N\in\N$ and positive increasing sequence $\{r_n\}$ with ${\displaystyle \lim_{n \to \infty}r_n = \infty}$ such that
\begin{equation}
 \limsup_{r \to \infty}\frac{\zeta_{\rho}^{2*}(x>r)}{\zeta_{\rho}(x >r)} < e^{-\delta}2^{N+1},\label{Prop5.1.4}
\end{equation}
and 
\begin{equation*}
\lim_{n \to \infty}\frac{\widetilde\eta^{N*}(|x|>r_n)}{\rho(|x| >r_n)} =0.
\end{equation*}
Define $I_j(r)$ and $J_j(r)$ for $j=1,2$ as
\begin{equation*}
I_j(r):=e^{-j\delta}\sum_{k=1}^{N}\frac{(j\delta)^k\widetilde\eta^{k*}(|x|>r)}{k!}\quad\mbox{and}\quad J_j(r):=e^{-j\delta}\sum_{k=N+1}^{\infty}\frac{(j\delta)^k\widetilde\eta^{k*}(|x|>r)}{k!}.
\end{equation*}
By Lemma \ref{lem5.1} (ii), there are $s_1 >0$ and $c_4 >0$ such that, for $1 \le k \le N$,
\begin{equation*}
\widetilde\eta^{k*}(|x|>r+s_1) \le c_4 \widetilde\eta^{N*}(|x|>r) \q \mbox{ for } r \ge 0.
\end{equation*}
Since $\rho \in {\mathcal O}{\mathcal S}$, we find from Lemma \ref{lem5.5} (i)  that there is $c_5 >0$ such that
 \begin{equation*}
\rho(|x| >r) \le c_5\rho(|x|>r+s_1)\q \mbox{ for } r \ge 0.
\end{equation*}
Thus we have, for $1 \le k \le N$, 
\begin{equation*}
\limsup_{n \to \infty}\frac{\widetilde\eta^{k*}(|x|>r_n+s_1)}{\rho(|x|>r_n+s_1)}\le  c_4 c_5 \lim_{n \to \infty}\frac{\widetilde\eta^{N*}(|x|>r_n)}{\rho(|x| >r_n)} =0.
\end{equation*}
Here, as we have $I_1(r)+J_1(r)=\rho(|x|>r)$ and $I_1(r_n+s_1)/\rho(|x|>r_n+s_1)\to 0$ as $n\to\infty$, then
$$\lim_{n\to\infty}\frac{J_1(r_n+s_1)}{\rho(|x|>r_n+s_1)}=1.$$
 Hence we establish that
\begin{eqnarray*}
\limsup_{r \to \infty}\frac{\zeta_{\rho}^{2*}(x>r)}{\zeta_{\rho}(x >r)} &\ge& \limsup_{r \to \infty}\frac{\rho^{2*}(|x|>r)}{\rho(|x| >r)}\nonumber\\
 &\ge& \limsup_{n \to \infty}\frac{I_2(r_n+s_1)+J_2(r_n+s_1)}{\rho(|x| >r_n+s_1)}\nonumber\\
 &\ge&\liminf_{r \to \infty}\frac{J_2(r)}{J_1(r)}\ge e^{-\delta}2^{N+1}. \nonumber
\end{eqnarray*}
This  contradicts \eqref{Prop5.1.4}. Thus \eqref{Prop5.1.3} holds and it follows from \eqref{Prop5.1.1} that there is  $n\in\N$ such that \eqref{Prop5.1.2} holds. We find from Lemma \ref{lem5.3} (i) that  $\widetilde\eta^{n*}\in {\mathcal O}{\mathcal S}$. 

Assertion (ii) is clear from (i). Next we prove (iii).  Suppose that $\rho \in {\mathcal O}{\mathcal S}$. Define a compound Poisson distribution $\tau$ on $\R_+$ such that its L\'evy measure is $\delta \zeta$. Recall that $\delta=\eta(\R^d)$. By Lemma \ref{lem5.1} (i), we have
\begin{eqnarray*}
 \zeta_{\rho}(x >r) = e^{-\delta}\sum_{k=1}^{\infty}\frac{\delta^k\widetilde\eta^{k*}(|x|>r)}{k!} \le  \tau(x >r) \mbox{ for } r\ge 0. 
\end{eqnarray*}
Note that $ \zeta_{\rho}(x >r)\ge c^{-1}\zeta(x >r)$ with $c:=
e^{\delta}\delta^{-1}$. By using \eqref{eq:5.8}, we see from Lemma \ref{lem5.2} that, for $ r\ge 0,$ 
\begin{eqnarray*}
\tau(x >r) &=& e^{ -\delta}\sum_{k=1}^{\infty}\frac{\delta^{k}}{k!}\zeta^{k*}(x>r)\\
             &\le& e^{ -\delta}\sum_{k=1}^{\infty}\frac{(c\delta)^k}{k!}\zeta_{\rho}^{k*}(x>r) \nonumber\\
             &\le& c_1 e^{ -\delta}\sum_{k=1}^{\infty}\frac{(c\delta)^{k}}{k!}(\ell^*(\zeta_{\rho})-1 +\ep)^{k}\zeta_{\rho}(x>r).\nonumber\
\end{eqnarray*}
Thus we find that 
\[
\tau(x >r) \asymp \zeta_{\rho}(x >r) =  \rho(|x| >r)
\]
and, by Lemma \ref{lem5.3} (i),  $\tau\in {\mathcal O}{\mathcal S}$ on $\R_+.$ Use (i) for the compound Poisson distribution $\tau$ on $\R_+$. There is a positive integer $m$ such that $\zeta^{m*}\in {\mathcal O}{\mathcal S}$ on $\R_+$ and \eqref{eq:3.14} holds. By Lemma \ref{lem5.3}, the converse assertion is trivial. \qed

\medskip

Now we prove Theorem \ref{th3.1}.

\medskip

\noindent
{\it Proof of Theorem \ref{th3.1}:} The proofs of (i), (ii), and (iii) are clear from Lemma \ref{lem5.6} and Proposition \ref{prop5.1}.  Assertion (iv) is proved in the same way of  Corollary 1.1 (ii) of \cite{sw05}. \qed

\bibliographystyle{alea2}
\bibliography{07-05}

\end{document}